\documentclass{amsart}
\usepackage{amssymb}
\usepackage{amsbsy}
\usepackage{amscd}
\usepackage{amsmath}
\usepackage{amsthm}

\theoremstyle{plain}
 \newtheorem{thm}{Theorem}[section]
 \newtheorem{lem}[thm]{Lemma}
 \newtheorem{prop}[thm]{Proposition}
 \newtheorem{cor}[thm]{Corollary}
 
\theoremstyle{definition}
 \newtheorem{defn}{Definition}[section]
\theoremstyle{remark}
 \newtheorem{rem}{Remark}[section]
 
 \newtheorem{claim}{Claim}[section]
\def\Bbb{\mathbb}
\def\frak{\mathfrak}
\def\cal{\mathcal}

\newcommand{\Ext}{\operatorname{Ext}}
\newcommand{\Hom}{\operatorname{Hom}}

\newcommand{\im}{\operatorname{im}}

\newcommand{\rk}{\operatorname{rk}}

\newcommand{\NS}{\operatorname{NS}}
\newcommand{\coker}{\operatorname{coker}}
\newcommand{\Pic}{\operatorname{Pic}}
\newcommand{\ch}{\operatorname{ch}}
\newcommand{\td}{\operatorname{td}}

\newcommand{\Quot}{\operatorname{Quot}}

\newcommand{\WIT}{\operatorname{WIT}}

\newcommand{\PD}{\operatorname{PD}}

\font\b=cmr10 scaled \magstep5
\def\bigzerou{\smash{\lower1.7ex\hbox{\b 0}}}

\numberwithin{equation}{section}
\begin{document}

\title{
Singularities on the 2-dimensional moduli spaces of stable sheaves on K3
surfaces}

\author{Nobuaki Onishi and K\={o}ta Yoshioka}
 
\address{
Department of mathematics, Faculty of Science, Kobe University,
Kobe, 657, Japan}

\email{onishi@math.kobe-u.ac.jp, yoshioka@math.kobe-u.ac.jp}

 \maketitle

\section{Introduction}

Let $X$ be a K3 surface over ${\Bbb C}$.
Mukai introduced a lattice structure $\langle \quad,\quad \rangle$
on 
$H^*(X,{\Bbb Z}):=\bigoplus_i H^{2i}(X,{\Bbb Z})$
by 
\begin{equation}
\begin{split}
\langle x,y \rangle:=&-\int_X x^{\vee} \wedge y\\
=& \int_X(x_1 \wedge y_1-x_0 \wedge y_2-x_2 \wedge y_0),
\end{split} 
\end{equation}
where $x_i \in H^{2i}(X,{\Bbb Z})$ (resp. $y_i \in H^{2i}(X,{\Bbb Z})$)
is the $2i$-th component of $x$ (resp. $y$)
and $x^{\vee}=x_0-x_1+x_2$.
It is now called the Mukai lattice.
For a coherent sheaf $E$ on $X$,
we can attach an element of $H^*(X,{\Bbb Z})$
called the Mukai vector 
\begin{equation}
v(E):=\ch(E)\sqrt{\td_X}
=\ch(E)(1+\rho_X),
\end{equation}
where $\ch(E)$ is the Chern character of $E$, 
$\td_X$ the Todd class of $X$ and $\rho_X$ the fundamental 
cohomology class of $X$ ($\int_X \rho_X=1$).

\begin{defn}\cite{Y:10}
We fix an ample divisor $H$ on $X$ and an element 
$G \in K(X) \otimes{\Bbb Q}$ with $\rk G>0$.
\begin{enumerate}
\item
Let $E$ be a torsion free sheaf on $X$.
$E$ is $G$-twisted semi-stable (resp. stable) with respect to $H$, if
\begin{equation}\label{eq:twist}
\frac{\chi(G,F(nH))}{\rk(F)} \leq \frac{\chi(G,E(nH))}{\rk(E)},
n \gg 0
\end{equation}
for $0 \subsetneq F \subsetneq E$
(resp. the inequality is strict).
\item
For a $w \in H^*(X,{\Bbb Q})_{alg}:={\Bbb Q} \oplus \NS(X) \otimes {\Bbb Q}
\oplus {\Bbb Q}\rho_X$ with $\rk w>0$,
we define the $w$-twisted semi-stability as the $G$-twisted semi-stability,
where $G \in K(X) \otimes {\Bbb Q}$ satisfies $v(G)=w$.
\end{enumerate}
\end{defn} 

Matsuki and Wentworth \cite{M-W:1} constructed the moduli space 
of $w$-twisted semi-stable sheaves $E$ with $v(E)=v$.
We denote it by $\overline{M}_H^w(v)$.
If $w=v({\cal O}_X)$, then the $v({\cal O}_X)$-twisted semi-stability
is nothing but the usual Gieseker's semi-stability.
Hence we denote $\overline{M}_H^{v({\cal O}_X)}(v)$ by
$\overline{M}_H(v)$. 

Assume that $v$ is an isotropic Mukai vector.
In \cite{Abe:1}, Abe considered the singularities of $\overline{M}_H(v)$.
Replacing $\overline{M}_H(v)$ by $\overline{M}_H^v(v)$,
we shall generalize Abe's results: 

\begin{thm}
\begin{enumerate}
\item[(1)]
$\overline{M}_H^v(v)$ is normal.
\item[(2)]  
For a suitable choice of $\alpha$ with $|\langle \alpha^2 \rangle| \ll 1$,
there is a surjective morphism
$\phi_{\alpha}:\overline{M}_H^{v+\alpha}(v)=
M_H^{v+\alpha}(v) \to \overline{M}_H^v(v)$ 
which becomes a minimal resolution of the singularities.
\item [(3)]
Let $x$ be a point of $\overline{M}_H^{v}(v)$
corresponding to the $S$-equivalence class
$\bigoplus_{i=0}^n E_i^{\oplus a_i}$, where $E_i$, 
$0 \leq i \leq n$ are $v$-twisted stable sheaves.
Then the matrix $(-\langle v(E_i),v(E_j) \rangle)_{i,j=0}^n$ is of affine
type $\tilde{A}_n,\tilde{D}_n, \tilde{E}_n$.
Assume that $a_0=1$. Then 
the singularity of $\overline{M}_H^{v}(v)$ at $x$ is a 
rational double point of type $A_n,D_n,E_n$ according as the type of 
the matrix $(-\langle v(E_i),v(E_j) \rangle)_{i,j=1}^n$.
\end{enumerate}
\end{thm}
Moreover we shall show that the Weyl chamber of the corresponding finite Lie algebra appears as a parameter space of $\alpha$.

If the matrix is of type $\tilde{A}_n$,
then the assertion (1) is due to Abe \cite[Thm. 3.3]{Abe:1}. 
Moreover if $n=1,2$, then he showed the assertion (3).
The assertion (2) is also contained in \cite[Thm. 3.3]{Abe:1}.
The main point of the proof is due to Matsuki and Wentworth \cite{M-W:1}
or Ellingsrud and G\"{o}ttsche \cite{E-G:1}.
The $v$-twisted stability naturally appears in the study of
the Fourier-Mukai transforms.
In \cite{Y:10}, \cite{Y:12}, we studied the Fourier-Mukai transform and
showed that the Fourier-Mukai transform preserves the $v$-twisted 
semi-stability under suitable assumptions.
So the $v$-twisted semi-stability is important and
this is our original motivation to study the moduli space of $v$-twisted 
semi-stable sheaves.
Another motivation is the following:
For the GIT quotients related to the moduli spaces of vector bundles 
on curves with additional structures,
the wall crossing behaviors have been studied by several authors.
In particular, Thaddeus \cite{Ta:1} described the wall crossing behavior
as a sequence of blowing-ups and blowing-downs and 
used it to show the Verlinde formula.
For the rank two case, Ellingsrud and G\"{o}ttsche \cite{E-G:1} studied
the similar variation problem for the  
moduli space of stable sheaves on a K3 surface.
In this case, Mukai's elementary transformation appears.
For all these examples, the exceptional locus of the blowing-up
is irreducible.
So it is interesting to construct an example with a reducible exceptional
lucus, and a rational double point will be a simple and interesting 
example to consider.

Our main idea to study the exceptional locus is
the same as the one in \cite{Y:5} to study the
Brill-Noether locus of sheaves on K3 surfaces.
Let us explain the contents of this note.
In section \ref{sect:defn}, we show that the 
Donaldson's determinant line bundle on $\overline{M}_H^v(v)$ is ample.
The $w$-twisted stability depends on the choice of $w$.
Hence we introduce a suitable parameter space of $w$ and  
 introduce a chamber structure on this parameter space.
This chamber will become a Weyl chamber of a finite simple Lie algebra
in the next section.
By using a special kind of Fourier-Mukai transform
called reflection,
we compare the twisted stability for two $w$
(Proposition \ref{prop:FM}).
Section \ref{sect:resolution} is the main part of this note.
We first describe the exceptional locus of the resolution
as a Brill-Noether locus of sheaves, under the assumption
that $w=v+\alpha$, $|\langle \alpha^2 \rangle| \ll 1$ belongs to a 
special chamber.  
By using the Weyl group action on the parameter space,
we can give a set-theoretic description of the exceptional locus
for general cases.
Finally we shall prove that $\overline{M}_H^v(v)$ is normal.
In section \ref{sect:example}, we give some examples of singular
moduli spaces by using the surjectivity of the period map.

\section{Definitions}\label{sect:defn}

Let $L$ be a lattice (or a ${\Bbb Q}$-vector space with a bilinear form)
with a weight 2 Hodge structure:
$L \otimes {\Bbb C}=\bigoplus_{p+q=2}L^{p,q}$.
We set $L_{alg}:=L \cap L^{1,1}$.
The Mukai lattice $H^*(X,{\Bbb Z})$ has a Hodge structure:
\begin{equation}
\begin{split}
H^{2,0}(H^*(X,{\Bbb C}))&=H^{2,0}(X),\\
H^{1,1}(H^*(X,{\Bbb C}))&=H^{0,0}(X) \oplus H^{1,1}(X) \oplus H^{2,2}(X),\\
H^{0,2}(H^*(X,{\Bbb C}))&=H^{0,2}(X).
\end{split}
\end{equation}
Then $H^*(X,{\Bbb Z})_{alg}={\Bbb Z} \oplus \Pic(X) \oplus {\Bbb Z}\rho_X$.

\subsection{Twisted stability}

Let $G$ be an element of $K(X) \otimes {\Bbb Q}$ with $\rk G>0$.
We fix an ample divisor $H$ on $X$.
For a coherent sheaf $E$ on $X$,
we define the $G$-twisted rank, 
degree, and Euler characteristic of $E$ by
\begin{equation}
\begin{split}
\rk_{G}(E)&:=\rk(G^{\vee} \otimes E),\\
\deg_{G}(E)&:=(c_1(G^{\vee} \otimes E),H),\\
\chi_{G}(E)&:=\chi(G^{\vee} \otimes E).
\end{split}
\end{equation}
We shall rewrite the condition \eqref{eq:twist} on the twisted stability.
By the Riemann-Roch theorem, we get that
\begin{equation}\label{eq:twist2}
\begin{split}
\frac{\chi(G,E(nH))}{\rk (G)\rk(E)}-\frac{\chi(G,F(nH))}{\rk(G)\rk(F)}=&
n \left(\frac{\deg_G(E)}{\rk_G(E)}-\frac{\deg_G(F)}{\rk_G(F)} \right)+
\left(\frac{\chi_G(E)}{\rk_G(E)}-\frac{\chi_G(F)}{\rk_G(F)}\right)\\
=&n \left(\frac{(c_1(E),H)}{\rk(E)}-\frac{(c_1(F),H)}{\rk(F)} \right)+
\left(\frac{\chi(E)}{\rk(E)}-\frac{\chi(F)}{\rk(F)}\right)\\
&
+\left(\frac{c_1(E)}{\rk(E)}-\frac{c_1(F)}{\rk(F)},
\frac{c_1(G)}{\rk G} \right).
\end{split}
\end{equation}
Let $\varphi:\Pic(X) \otimes {\Bbb Q} \to H^{\perp}$ be the orthogonal
projection.
Then the twisted stability depends only on 
$\varphi(c_1(G)/\rk G) \in H^{\perp}$
and it is nothing but the twisted stability due to
Matsuki-Wentworth \cite{M-W:1}. 
\begin{defn}\label{defn:general}
A polarization $H$ is general with respect to $v$, 
if the following condition holds:
\begin{itemize}
\item[$(*)$]
for every $\mu$-semi-stable sheaf $E$ with $v(E)=v$,
if $F \subset E$ satisfies
$(c_1(F),H)/\rk F=(c_1(E),H)/\rk E$, then
$c_1(F)/\rk F=c_1(E)/\rk E$.
\end{itemize} 
\end{defn}
If $H$ is general with respect to $v$, then
the $w$-twisted semi-stability does not depend on the choice of $w$.  
The following theorem was proved in \cite{M-W:1}.
\begin{thm}\cite{M-W:1}\label{thm:M-W}
Let $w$ be an element 
of $H^*(X,{\Bbb Q})_{alg}$ such that $\rk w>0$. 
Then there is a coarse moduli scheme $\overline{M}_H^{w}(v)$
of $S$-equivalence classes of $w$-twisted semi-stable sheaves
$E$ with $v(E)=v$.
$\overline{M}_H^{w}(v)$ is a projective scheme.
\end{thm}

\begin{defn}
We denote the open subscheme of $\overline{M}_H^{w}(v)$
consisting of $w$-twisted stable sheaves
by $M_H^{w}(v)$.  
\end{defn}
If $w=v({\cal O}_X)$, then we denote $\overline{M}_H^w(v)$
(resp. ${M}_H^w(v)$) by
$\overline{M}_H(v)$ (resp. ${M}_H(v)$).

\begin{prop}\label{prop:exist}
$\overline{M}_H^w(v) \ne \emptyset$ if $\langle v^2 \rangle \geq -2$.
\end{prop}

\begin{proof}
We may assume that $v$ is primitive.
If $H$ is general with respect to $v$, then
\cite[Thm. 8.1]{Y:7} implies that $\overline{M}_H(v) \ne \emptyset$.
By the study of the 
chamber structure (cf. \cite[sect. 1, Prop. 4.2]{Y:10}),
we get our claim
\end{proof}

\subsection{Line bundles on $\overline{M}_H^{w}(v)$}
Throughout this note,
$v:=r+\xi+a \rho_X$, $\xi \in \Pic(X)$ is a primitive isotropic
Mukai vector with $r>0$.

We define a homomorphism which preserves the Hodge structure and the metric: 
\begin{equation}\label{eq:delta}
\begin{matrix}
\delta:H^2(X,{\Bbb Q})& \to & H^*(X,{\Bbb Q})\\
D & \mapsto & D+\frac{(D,\xi)}{r}\rho_X.
\end{matrix}
\end{equation}
We denote $\delta(D)$ by $\widehat{D}$.
Then we have an orthogonal decomposition:
\begin{equation}\label{eq:orthgonal}
H^*(X,{\Bbb Q})=
({\Bbb Q}v \oplus {\Bbb Q}\rho_X) \bot \delta(H^2(X,{\Bbb Q})).
\end{equation} 
Let $\theta_v^{\alpha}: v^{\perp} \to H^2(M_H^{v+\alpha}(v),{\Bbb Z})$
be the Mukai homomorphism defined by
\begin{equation}
 \theta_v^{\alpha}(x):=\frac{1}{\rho}
 \left[p_{M_H^{v+\alpha}(v)*}((\ch {\cal E})\sqrt{\td_X}
 x^{\vee})\right]_1,
\end{equation}
where ${\cal E}$ is a quasi-universal family of similitude 
$\rho$.
If $x \in (v^{\perp}/{\Bbb Z}v) \otimes {\Bbb Q}$, then we have a
${\Bbb Q}$-line bundle $L(x)$ on $M_H^{v+\alpha}(v)$ such that
$c_1(L(x))=\theta_v^{\alpha}(x)$.
For $\widehat{H}=H+\{(H,\xi)/r\} \rho_X$,
$L(r\widehat{H})$ is the Donaldson's determinant line bundle
and J. Li \cite{Li:1} showed that
canonically $L(r\widehat{H})$ 
extends to a line bundle on $\overline{M}_H^{v+\alpha}(v)$.
We also denote this extension by $L(r\widehat{H})$.
Then $L(r\widehat{H})$ is a nef and big line bundle 
and we have a contraction map
from the Gieseker moduli space to the Uhlenbeck moduli space.   
Hence $L(r\widehat{H})$ is important.

One of the reason we consider the $v$-twisted  stability
is the following proposition.
\begin{prop}
$L(r\widehat{H})$ is an ample line bundle on $\overline{M}_H^v(v)$.
\end{prop}

\begin{proof}
We recall the construction of $\overline{M}_H^v(v)$ in
\cite{Y:11}.
Let $E$ be a $v$-twisted stable sheaf with $v(E)=v$.
We set $N:=\chi(E,E(nH))$.
Let $Q:=\Quot_{E(-nH)^{\oplus N}/X/{\Bbb C}}^v$ be a quot-scheme
parametrizing all quotients
$E(-nH)^{\oplus N} \to F$ such that
$v(F)=v$ and 
${\cal O}_Q^{\oplus N} \otimes E(-nH) \to {\cal Q}$ the universal
quotient.
Let $Q^{ss}$ be an open subscheme of $Q$
consisting of $q \in Q$ such that
\begin{enumerate}
\item
 ${\cal Q}_q$ is $v$-twisted semi-stable,
\item
 $\Hom(E,E^{\oplus N}) \to \Hom(E,{\cal Q}_q(nH))$ is an isomorphism,
\item
 $\Ext^i(E,{\cal Q}_q(nH))=0$, $i>0$.
\end{enumerate}
Then we have an isomorphism
${\cal O}_{Q^{ss}}^{\oplus N} \to 
p_{Q^{ss}*}({\cal Q} \otimes p_X^*(E(-nH))^{\vee})$.
$Q$ has a natural action of $GL(N)$.
We set ${\cal L}_{m}:=\det p_{Q!}({\cal Q} \otimes p_X^*(E(-mH))^{\vee})$.
Since ${\cal Q}$ is $GL(N)$-linearized, 
${\cal L}_m$ is also $GL(N)$-linearized.
By the construction of $Q$,
 ${\cal L}_{n+m}$, $m \gg 0$ gives an embedding of $Q$ to a
Grassmann variety.
Thus ${\cal L}_{n+m}$, $m \gg 0$ is ample.
Let $T:=\det ({\cal O}_{Q}^{\oplus N})$ be the $GL(N)$-linearized line bundle
induced by the standard action of $GL(N)$ on
${\cal O}_{Q}^{\oplus N}$.
The center ${\Bbb C}^{\times} \subset GL(N)$ acts trivially on $Q$ and
the action on ${\cal L}_m$ is the multiplication by
$\chi(E,E({\cal O}_X(-mH)))$-th power of constant.
By a simple calculation, we see that 
$\chi(E,E({\cal O}_X(-mH)))=r^2m^2(H^2)/2$.
Hence ${\cal L}:=
{\cal L}_{n+m}^{\otimes m^2} \otimes {\cal L}_{n}^{\otimes -(n+m)^2}$
and ${\cal L}':={\cal L}_{n+m}^{\otimes m^2} \otimes T^{\otimes -(n+m)^2}$
have $PGL(N)$-linearizations.
By the construction of the moduli space,
$\overline{M}_H^v(v)$ is described as a GIT quotient
$Q^{ss} \to \overline{M}_H^v(v)$, where $n \gg 0$ and
${\cal L}'$ is the linearization.
Since ${\cal L}_{|Q^{ss}}={\cal L}'_{|Q^{ss}}$ as $PGL(N)$-line bundles,
${\cal L}_{|Q^{ss}}$ descends to an ample line bundle on
$\overline{M}_H^v(v)$. 
We note that 
\begin{equation}
{\cal L}={\cal L}_{n+m}^{\otimes m^2} \otimes {\cal L}_{n}^{\otimes -(n+m)^2}  
=\det p_{Q!}({\cal Q} \otimes p_X^*(L)^{\vee}), 
\end{equation}
where
$L=m^2 E(-(n+m)H)-(n+m)^2 E(-nH) \in K(X)$.
Since $\det p_{Q^{ss}!}({\cal Q} \otimes p_X^*(E)^{\vee})={\cal O}_{Q^{ss}}$,
\begin{equation}
\det p_{Q^{ss}!}({\cal Q} \otimes p_X^*(L)^{\vee})=
\det p_{Q^{ss}!}({\cal Q} \otimes p_X^*({L'})^{\vee}),
\end{equation}
where $L'=L-(m^2-(n+m)^2)E$.
Since $v(L')=rmn(m+n) \widehat{H}$,
we get our claim.
\end{proof}

\begin{cor}
\begin{enumerate}
\item [(1)]
 If $\xi \in {\Bbb Q}H$, then $\overline{M}_H(v)=\overline{M}_H^v(v)$.
Hence $L(r\widehat{H})$ is an ample line bundle on $\overline{M}_H(v)$.
\item
[(2)]
 Let ${\cal M}_n$ (resp. $\overline{\cal M}_n$ )
be the moduli space of polarized (resp. quasi-polarized) 
K3 surfaces $(X,H)$ with $(H^2)=2n$.
We set $v:=r+dH+a \rho_X$, $d^2(H^2)=2ra$.
Then we have a morphism of the moduli spaces 
${\cal M}_n \to \overline{\cal M}_{n'}$:
$(X,H) \mapsto (\overline{M}_H(v),L(\widehat{H}))$
where $n=ra/d^2$ and $n'$ is determined by the primitive class
in ${\Bbb Q}L(\widehat{H}) \cap \Pic(\overline{M}_H(v))$.
\end{enumerate}
\end{cor}
In particular, if 
$v:=r+H+a \rho_X$ satisfies $\gcd(r,a)=1$, then
$M_H(v)$ is compact and
$\widehat{H}=H+2a \rho_X$ gives a canonical primitive polarization of 
$M_H(v)=M_H^v(v)$.   

\begin{rem}
If $\gcd(r,d)=1$, then $M_H(v)$ consists of
$\mu$-stable locally free sheaves for a general $X$.
For a special $X$,
$M_H(v)$ may consist of properly $\mu$-semi-stable
sheaves. Indeed
let $X \to {\Bbb P}^1$ be an elliptic K3 surface with a 
section $\sigma$.
We set $H:=\sigma+3f$, where $f$ is a fiber of $\pi$. 
If $\Pic(X)={\Bbb Z} \sigma \oplus {\Bbb Z} f$,
then $H$ is an ample divisor with $(H^2)=4$.
We set $v=2+H+\rho_X$. Then $\langle v^2 \rangle=0$
and every member of $M_H(v)$ is given by
\begin{equation}
E:=\ker({\cal O}_X(\sigma+f) \oplus {\cal O}_X(2f) 
\overset{ev}{\to} {\Bbb C}_s),
\end{equation}
$s \in X$. 
\end{rem}
 
\vspace{1pc}

For the Mukai homomorphism, Mukai \cite{Mu:4} showed the following.
\begin{thm}[Mukai] 
Assume that $M_H^{v+\alpha}(v)$ is compact. Then  
$\theta_v^{\alpha}$ is surjective and 
the kernel is ${\Bbb Z}v$.
Moreover 
$\theta_v^{\alpha}:v^{\perp}/{\Bbb Z}v \to H^2(M_H^{v+\alpha}(v),{\Bbb Z})$
is a Hodge isometry.
\end{thm}

By \eqref{eq:delta} and \eqref{eq:orthgonal},
we have a sequence of Hodge isometries:
\begin{equation}\label{eq:hodge1}
H^2(X,{\Bbb Q}) \to \delta(H^2(X,{\Bbb Q})) \to
(v^{\perp}/{\Bbb Z}v) \otimes_{\Bbb Z} {\Bbb Q} \to
H^2(M_H^{v+\alpha}(v),{\Bbb Q}).
\end{equation}
Then since $\widehat{H}=\delta(H) \in v^{\perp}$,
we have an isometry
\begin{equation}\label{eq:hodge2}
\delta(H^{\perp})_{alg} \to 
((\widehat{H}^{\perp} \cap v^{\perp})/{\Bbb Z}v)_{alg} \otimes {\Bbb Q}.
\end{equation}
In particular $((\widehat{H}^{\perp} \cap v^{\perp})/{\Bbb Z}v)_{alg}$
is negative definite.

\subsection{Chamber structure}

We shall study the dependence of $\overline{M}_H^w(v)$ on $w$.
By \eqref{eq:twist2},
we may assume that $w=v+\alpha, \alpha \in \delta(H^{\perp})_{alg}$.
Let $u$ be a Mukai vector such that $0<\rk u<\rk v$, 
$\langle v,u \rangle \leq 0$, $\langle u^2 \rangle=-2$ and
$\langle u,\widehat{H} \rangle=0$.
We define a wall $W_{u} \subset 
\delta(H^{\perp})_{alg} \otimes_{\Bbb Q} {\Bbb R}$ 
with respect to $v$
by 
\begin{equation}
W_u:=\{\alpha \in \delta(H^{\perp})_{alg} \otimes {\Bbb R}|\;
\langle v+\alpha,u \rangle=0 \}.
\end{equation}

For a properly $v+\alpha$-twisted semi-stable sheaf $E$
with $v(E)=v$,
we consider the Jordan-H\"{o}lder filtration
\begin{equation}
 0 \subset F_1 \subset F_2 \subset \dots \subset F_s=E
\end{equation}
of $E$ with respect to the $v+\alpha$-twisted stability, that is, 
$E_i:=F_i/F_{i-1}$ is a $v+\alpha$-twisted stable sheaf with
\begin{equation}
\begin{split}
&(c_1(E_i),H)/\rk E_i=(c_1(E),H)/\rk E=(\xi,H)/r,\\
&\langle v+\alpha,v(E_i) \rangle/\rk E_i=\langle v+\alpha,v \rangle/ \rk E.
\end{split}
\end{equation}
We set $u_i:=v(E_i)$.
Then we see that
$\langle u_i, \widehat{H} \rangle=
(c_1(E_i),H)-(\xi,H)r_i/r=0$ and
$\langle v+\alpha,u_i \rangle=0$.

\begin{lem}
$\langle v, u_i \rangle \leq 0$ and
$\langle u_i^2 \rangle=-2$ for some $i$.
\end{lem}

\begin{proof}
Since $u_i/\rk u_i-v/\rk v \in \widehat{H}^{\perp} \cap \rho_X^{\perp}$
and $(\widehat{H}^{\perp} \cap \rho_X^{\perp})_{alg}$ 
is negative semi-definite,
$\langle (u_i/\rk u_i-v/\rk v)^2 \rangle \leq 0$.
Then $\langle u_i^2 \rangle \leq 2 \langle u_i,v \rangle (\rk u_i/\rk v)$.
Since $\sum_j \langle u_j,v \rangle=\langle v,v \rangle=0$,
we get $\langle u_i,v \rangle \leq 0$ for some $i$.
In particular $\langle u_i^2 \rangle<0$ provided that 
$\langle u_i,v \rangle < 0$.
If $\langle u_i,v \rangle=0$, then
$u_i \in  \widehat{H}^{\perp} \cap v^{\perp}$.
Since $\rk u_i<\rk v$, we get $u_i \not \in {\Bbb Z}v$. Then
\eqref{eq:hodge2} implies that
$\langle u_i^2 \rangle<0$.
Since $\langle u_i^2 \rangle \geq -2$, we conclude that
$\langle u_i^2 \rangle=-2$.
\end{proof}
Therefore $\alpha \in W_{u_i}$.
We set 
\begin{equation}
 {\cal U}:=\left\{u \in H^*(X,{\Bbb Z})_{alg}\left|
\begin{split}
&\langle u^2 \rangle=-2, \langle v,u \rangle \leq 0,
\langle \widehat{H},u \rangle=0,\\
&0< \rk u <\rk v 
\end{split}
\right. \right\}.
\end{equation}
For a fixed $v$ and $H$, ${\cal U}$ is a  finite set.

\begin{lem}
If $\alpha$ does not lie on any wall $W_u$, $u \in {\cal U}$,
then $\overline{M}_H^{v+\alpha}(v)=M_H^{v+\alpha}(v)$.
In particular, $\overline{M}_H^{v+\alpha}(v)$ is a K3 surface.
\end{lem} 

\begin{defn}
Let ${\cal C}$ be a connected component of 
$\delta(H^{\perp})_{alg} \otimes_{\Bbb Q} {\Bbb R} 
\setminus \cup_{u \in {\cal U}} W_u$.
We call ${\cal C}$ a chamber.
\end{defn}
As is proved in
\cite{M-W:1}, we get
\begin{prop}
 The $v+\alpha$-twisted stability does not depend on the choice of 
$\alpha \in {\cal C}$.
If $\beta$ belongs to the closure of ${\cal C}$,
then we have a morphism $\overline{M}_H^{v+\alpha}(v) \to 
\overline{M}_H^{v+\beta}(v)$ for $\alpha \in {\cal C}$.
In particular, we have a morphism $\phi_{\alpha}:\overline{M}_H^{v+\alpha}(v)
\to \overline{M}_H^{v}(v)$ for $|\langle \alpha^2 \rangle| \ll 1$.
\end{prop}
Let $T \subset v^{\perp}$ be a sufficiently small neighborhood of $0$.
Then $W_u$ 
intersects $T$ if and only if $\langle v,u \rangle=0$.
Since we are interested in the neighborhood of $v$,
we may assume that the defining equation of a wall
$W_u$ belongs to the subset
\begin{equation}
{\cal U}':=\{u \in {\cal U}|\;
\langle v,u \rangle=0 \}.
\end{equation}

By the same argument as above, we get the following.
\begin{lem}\label{lem:JHF}  
Let $E$ be a properly $v$-twisted semi-stable sheaf with $v(E)=v$ and
\begin{equation}
0 \subset F_1 \subset F_2 \subset \dots \subset F_s=E
\end{equation}
the Jordan-H\"{o}lder filtration of $E$ with respect to
the $v$-twisted stability. 
Then $\langle v(F_i/F_{i-1})^2 \rangle=-2$.
\end{lem}

\subsection{Reflection}
For an $\alpha \in \delta(H^{\perp})_{alg}$ with 
$|\langle \alpha^2 \rangle| \ll 1$, 
let $F$ be a $v+\alpha$-twisted stable torsion free sheaf such that
\begin{enumerate}
\item
$\langle v(F)^2 \rangle=-2$,
\item
$\langle v(F),\widehat{H} \rangle/\rk F
=(c_1(F),H)/\rk F-(\xi,H)/r=0$ and 
\item
$\langle v,v(F) \rangle=\langle \alpha, v(F)\rangle=0$.
\end{enumerate}
By (i), $F$ is a rigid torsion free sheaf, and hence $F$ is locally free.

Let ${\cal E}$ be a coherent sheaf on $X \times X$ which is
defined by an exact sequence 
\begin{equation}\label{eq:kernel}
0 \to {\cal E} \to p_1^*(F^{\vee}) \otimes p_2^*(F) \overset{ev}{\to}
 {\cal O}_{\Delta} \to 0, 
\end{equation}
where $p_i:X \times X \to X$, $i=1,2$ are projections.
We consider the Fourier-Mukai transform induced by ${\cal E}$:
\begin{equation}
\begin{matrix}
{\cal F}_{\cal E}:&{\bf D}(X) & \to & {\bf D}(X)\\
& x & \mapsto & {\bf R}p_{2*}(p_1^*(x) \otimes {\cal E}),
\end{matrix}
\end{equation}  
where ${\bf D}(X)$ is the bounded derived category of $X$.
Up to shift, the inverse of ${\cal F}_{\cal E}$ is given by 
\begin{equation}
\begin{matrix}
\widehat{\cal F}_{\cal E}:&{\bf D}(X) & \to & {\bf D}(X)\\
& y & \mapsto & {\bf R}\Hom_{p_{1}}( {\cal E}, p_2^*(y)).
\end{matrix}
\end{equation}  

\begin{defn}
Let $E$ be a coherent sheaf on $X$.
\begin{enumerate}
\item
We denote the $i$-th cohomology sheaf of
${\cal F}_{\cal E}(E)$ (resp. $\widehat{\cal F}_{\cal E}(E)$)
by ${\cal F}_{\cal E}^i(E)$ (resp. $\widehat{\cal F}_{\cal E}^i(E)$).  
\item
$E$ satisfies $\WIT_i$ with respect to ${\cal F}_{\cal E}$
(resp. $\widehat{\cal F}_{\cal E}$),
if ${\cal F}_{\cal E}^j(E)=0$ (resp. $\widehat{\cal F}_{\cal E}^j(E)=0$)
for $j \ne i$.
\end{enumerate}
\end{defn}

The Fourier-Mukai transform ${\cal F}_{\cal E}$ 
induces an isometry of the Mukai lattice 
${\cal F}_{\cal E}:H^*(X,{\Bbb Z}) \to H^*(X,{\Bbb Z})$.
Let $R_{v(F)}:H^*(X,{\Bbb Z}) \to H^*(X,{\Bbb Z})$ be the reflection
defined by the $(-2)$-vector $v(F)$:
\begin{equation}
 R_{v(F)}(u)=u+\langle u,v(F) \rangle v(F), u \in H^*(X,{\Bbb Z}).
\end{equation}
Then we see that ${\cal F}_{\cal E}=-R_{v(F)}$.
Thus the Fourier-Mukai transform ${\cal F}_{\cal E}$ is the geometric
realization of the reflection $R_{v(F)}$.

\begin{lem}\label{lem:vanish1}
Let $G$ be a $v+\alpha$-twisted semi-stable sheaf such that
$\deg(G)/\rk G=\deg(F)/\rk F$ and
$\chi(E+A,G) \geq 0$, where $E, A \in K(X) \otimes {\Bbb Q}$
satisfy $v(E)=v$,
$v(A)=\alpha$.
Then 
\begin{equation}
\Ext^2({\cal E}_{|\{x \} \times X},G)=0
\end{equation}
 for all $x \in X$.
\end{lem}

\begin{proof}
Assume that there is a non-zero homomorphism
$\varphi:G \to{\cal E}_{|\{x \} \times X}$.
Then we have a non-zero homomorphism $\psi:G \to F$.
Since $G$ is a $v+\alpha$-twisted semi-stable sheaf with
$\deg(G)/\rk G=\deg(F)/\rk F$ and
$F$ is $v+\alpha$-twisted stable,
we get that $0 \leq \chi(E+A,G)/\rk G \leq \chi(E+A,F)/\rk F=0$.
Hence
$\chi(E+A,G)=0$ and $\psi$ is surjective. 
Thus $\im \varphi$ contains $F$. 
On the other hand,
by the construction of ${\cal E}_{|\{x \} \times X}$,
${\cal E}_{|\{x \} \times X}$ does not contain $F$.
Therefore $\Hom(G,{\cal E}_{|\{x \} \times X})=0$.
By the Serre duality, we get
$\Ext^2({\cal E}_{|\{x \} \times X},G)=0$.
\end{proof}

\begin{lem}\label{lem:ev}
Let $E,A \in K(X) \otimes {\Bbb Q}$ be as in Lemma \ref{lem:vanish1}.
Let $G$ be a $v+\alpha$-twisted semi-stable sheaf 
such that $\deg(G)/\rk G=\deg(F)/\rk F$ and
$\chi(E+A,G)=0$.
Then the evaluation map
$\phi:\Hom(F,G) \otimes F \to G$ is injective and
$\coker \phi$ is a $v+\alpha$-twisted semi-stable sheaf.
\end{lem}

\begin{proof}
By the $v+\alpha$-twisted semi-stability of $E$ and $F$,
we see that $\deg(\im \phi )/\rk (\im \phi)=\deg(F)/\rk F$ and
$\chi(E+A,\im \phi)=0$.
Hence we get $\deg(\ker \phi)/\rk (\ker \phi)=\deg(F)/\rk F$ and
$\chi(E+A,\ker \phi)=0$.
Assume that $\ker \phi \ne 0$.
By the $v+\alpha$-twisted semi-stability of $\Hom(F,G) \otimes F$,
$\ker \phi$ is $v+\alpha$-twisted semi-stable.
Then we see that $\ker \phi \cong F^{\oplus k}$,
which implies that $\Hom(F,\ker \phi) \ne 0$.
On the other hand, $\phi$ induces an isomorphism
$\Hom(F,G) \otimes \Hom(F,F) \to \Hom(F,G)$.
Hence we have $\Hom(F,\ker \phi)=0$, which is a contradiction.
Therefore $\ker \phi=0$.
\end{proof}

\begin{prop}\label{prop:FM}
We set $\alpha^{\pm}:=\pm \epsilon v(F)+\alpha$, where
$0<\epsilon \ll 1$.
\begin{enumerate}
\item[(1)]
Let $E$ be a $v+\alpha^-$-twisted semi-stable sheaf with $v(E)=v$.
Then $\WIT_1$ holds for $E$ with respect to ${\cal F}_{\cal E}$
and ${\cal F}_{\cal E}^1(E)$ is a $v+\alpha^+$-twisted semi-stable
sheaf.
\item[(2)]
Conversely, for a $v+\alpha^+$-twisted semi-stable sheaf $E$ with
$v(E)=v$, $\WIT_1$ holds with respect to $\widehat{\cal F}_{\cal E}$
and $\widehat{\cal F}_{\cal E}^1(E)$ is a $v+\alpha^-$-twisted semi-stable
sheaf.
\item[(3)]
Moreover ${\cal F}_{\cal E}$ preserves the $S$-equivalence classes.
Hence we have an isomorphism
\begin{equation}
\overline{M}_H^{v+\alpha^-}(v) \to
\overline{M}_H^{v+\alpha^+}(v).
\end{equation}
\end{enumerate}
\end{prop}

\begin{proof}
We take an element $A \in K(X) \otimes {\Bbb Q}$ such that
$v(A)=\alpha$. We note that $F$ is 
$(v \pm \epsilon v(F)+\alpha)$-twisted stable
for $0 \leq \epsilon \ll 1$.
We first prove (1).
We note that $E$ is $v+\alpha^-$-twisted semi-stable.
By the definition of ${\cal E}$, we get an exact sequence
\begin{equation}
\begin{CD}
0 @>>> p_{2*}({\cal E} \otimes p_1^*(E)) @>>>
\Hom(F,E) \otimes F @>>> E \\
@>>> R^1 p_{2*}({\cal E} \otimes p_1^*(E)) @>>>
\Ext^1(F,E) \otimes F @>>> 0\\
@>>> R^2 p_{2*}({\cal E} \otimes p_1^*(E)) @>>>
\Ext^2(F,E) \otimes F @>>> 0.
\end{CD}
\end{equation}
Since $\deg F/\rk F=\deg E/\rk E$ and
$\chi(E-\epsilon F+A,F)/\rk F=-2\epsilon/\rk F<0
=\chi(E-\epsilon F+A,E)/\rk E$, the $v+\alpha^-$-twisted 
semi-stability of $E$ and $F$ imply that 
$\Ext^2(F,E)=\Hom(E,F)^{\vee}=0$. 
Thus $R^2 p_{2*}({\cal E} \otimes p_1^*(E))=0$. 
Since $E$ is $v+\alpha$-twisted semi-stable,
Lemma \ref{lem:ev} implies that
 $\Hom(F,E) \otimes F \to E$ is injective, and hence
$p_{2*}({\cal E} \otimes p_1^*(E))=0$.
Therefore $\WIT_1$ holds for $E$ and ${\cal F}_{\cal E}^1(E)$ is a 
$v+\alpha$-twisted
semi-stable sheaf with $v({\cal F}_{\cal E}^1(E))=R_{v(F)}(v)=v$.
Assume that ${\cal F}_{\cal E}^1(E)$ is not $v+\alpha^+$-twisted 
semi-stable.
Then there is an exact sequence
\begin{equation}
 0 \to G_1 \to {\cal F}_{\cal E}^1(E) \to G_2 \to 0
\end{equation}
such that $G_1$ is a $v+\alpha$-twisted semi-stable sheaf 
with $\deg_E(G_1)=\chi(E+A,G_1)=0$ and
$G_2$ is a $v+\alpha^+$-twisted stable sheaf with 
$\chi(E+\epsilon F+A,G_2)<0$.
By Lemma \ref{lem:vanish1}, we get 
$\widehat{\cal F}_{\cal E}^2(G_1)=0$.
Since $\chi(E+\epsilon F+A,G_2)<0$ and $\chi(E+\epsilon F+A,F)=2\epsilon>0$,
$\widehat{\cal F}_{\cal E}^0(G_2)=\Hom(F,G_2) \otimes F=0$.
Therefore $\WIT_1$ holds for $G_1$, $G_2$ and we get an exact sequence 
\begin{equation}
 0 \to \widehat{\cal F}_{\cal E}^1(G_1) \to E \to  
 \widehat{\cal F}_{\cal E}^1(G_2) \to 0.
\end{equation}
Since $\chi(F,\widehat{\cal F}_{\cal E}^1(G_2))=-\chi(F,G_2)>0$,
we get a contradiction.

(2) Conversely, let $E$ be a $v+\alpha^+$-twisted semi-stable sheaf with
$v(E)=v$. Then we have an exact sequence
\begin{equation}
\begin{CD}
0 @>>> \Hom_{p_{1}}({\cal O}_{\Delta}, p_2^*(E)) @>>>
\Hom(F,E) \otimes F @>>> \Hom_{p_{1}}({\cal E}, p_2^*(E)) \\
@>>> \Ext^1_{p_{1}}({\cal O}_{\Delta}, p_2^*(E)) @>>>
\Ext^1(F,E) \otimes F @>>> \Ext^1_{p_{1}}({\cal E}, p_2^*(E))\\
@>>> \Ext^2_{p_{1}}({\cal O}_{\Delta}, p_2^*(E)) @>>>
\Ext^2(F,E) \otimes F @>>> \Ext^2_{p_{1}}({\cal E}, p_2^*(E)).
\end{CD}
\end{equation}
By Lemma \ref{lem:vanish1},
$\widehat{\cal F}_{\cal E}^2(E)=\Ext^2_{p_{1}}({\cal E}, p_2^*(E))=0$.
It is easy to see that 
\begin{equation}
\begin{split}
&\Hom_{p_{1}}({\cal O}_{\Delta}, p_2^*(E))=0,\\
&\Ext^1_{p_{1}}({\cal O}_{\Delta}, p_2^*(E))=0,\\
&\Ext^2_{p_{1}}({\cal O}_{\Delta}, p_2^*(E))=E.
\end{split}
\end{equation}
Since $\chi(E+\epsilon F+A,E)=0<\chi(E+\epsilon F+A,F)$,
the $v+\alpha$-twisted semi-stability of $E$ and $F$ imply that
$\widehat{\cal F}_{\cal E}^0(E)=\Hom(F,E) \otimes F=0$.
Therefore $\WIT_1$ holds with respect to
$\widehat{\cal F}_{\cal E}$ and $\widehat{\cal F}_{\cal E}^1(E)$
is $v+\alpha$-twisted semi-stable.
Assume that $\widehat{\cal F}_{\cal E}^1(E)$ is not 
$v+\alpha^-$-twisted semi-stable.
Then there is an exact sequence
\begin{equation}
 0 \to G_1 \to \widehat{\cal F}_{\cal E}^1(E) \to G_2 \to 0
\end{equation}
such that $G_1$ is a $v+\alpha^-$-twisted stable sheaf with 
$\deg_E G_1=0$,
$\chi(E-\epsilon F+A,G_1)>0$, $\chi(E+A,G_1)=0$ and
$G_2$ is a $v+\alpha$-twisted semi-stable sheaf. 
Since $\chi(E-\epsilon F+A,G_1)/\rk G_1>0>\chi(E-\epsilon F+A,F)/\rk F$, 
$\Ext^2(F,G_1)=\Hom(G_1,F)^{\vee}=0$. 
Thus $R^2{p_{2*}}({\cal E} \otimes p_1^*(G_1))=0$. 
Since $G_2$ is a $v+\alpha$-twisted semi-stable sheaf with 
$\deg_E G_2=\chi(E+A,G_2)=0$, Lemma \ref{lem:ev} implies that
$\Hom(F,G_2) \otimes F \to G_2$ is injective, and hence 
$p_{2*}({\cal E} \otimes p_1^*(G_2))=0$.
Therefore $\WIT_1$ holds for $G_1$ and $G_2$ with respect to
${\cal F}_{\cal E}$.
Since $\chi(E+\epsilon F+A,{\cal F}_{\cal E}^1(G_1))>0$,
we get a contradiction.   

The last claim (3) will easily follow from the above arguments. 
We omit the proof.
\end{proof}

The following is proved in \cite{Y:5}.
\begin{prop}\label{prop:comm}
Keep notation as above.
Assume that $v+\alpha^-$ does not lie on walls.
Then 
${\cal F}_{\cal E}$ induces an isometry $R_{v(F)}:v^{\perp} \to v^{\perp}$
and the following diagram is commutative.
\begin{equation}\label{eq:diagram}
 \begin{CD}
  v^{\perp} @>{R_{v(F)}}>> v^{\perp}\\
  @V{\theta_v^{\alpha^-}}VV @VV{\theta_v^{\alpha^+}}V\\
  H^2(M_{H}^{v+\alpha^-}(v),{\Bbb Z}) @= 
  H^2(M_{H}^{v+\alpha^+}(v),{\Bbb Z})
 \end{CD}
\end{equation}
\end{prop}

\begin{rem}\label{rem:comm}
If $\alpha$ belongs to exactly one wall $W_u$, $u \in {\cal U}$, then
there is a $v+\alpha$-twisted stable sheaf $F$ with $v(F)=u$. 
So we can apply Propositions \ref{prop:FM} and
\ref{prop:comm} to this $F$.
\end{rem}

\section{Resolution of the singularities of $\overline{M}_H^{v}(v)$}
\label{sect:resolution}

\subsection{Exceptional locus of the resolution}

Assume that there is a point $x$ of $\overline{M}_H^{v}(v)$ 
representing a properly $v$-twisted semi-stable sheaf.
Assume that $x$ is represented by an $S$-equivalence class
$\bigoplus_{i=0}^{n} E_i^{\oplus a_i}$, where $E_i$ is a 
$v$-twisted stable sheaf such that
$\langle v(E_i),\widehat{H} \rangle=\langle v(E_i),v \rangle=0$
and $E_i \ne E_j$ for $i \ne j$.
We set $v_i:=v(E_i)$. 

\begin{lem}\label{lem:lattice}
${\Bbb Z}v_0+{\Bbb Z}v_1+\dots+{\Bbb Z}v_n$ 
is a negative semi-definite lattice of
affine type $\tilde{A}_n,\tilde{D}_n,\tilde{E}_n$.
More precisely, $(-\langle v_i,v_j \rangle)_{i,j=0}^n$ is the Cartan
matrix of the affine Lie algebra $\tilde{A}_n,\tilde{D}_n,\tilde{E}_n$.
In particular, $v_0,v_1,\dots,v_n$ are linearly independent.
\end{lem}

\begin{proof}
We note that 
\begin{equation}
\begin{cases}
\langle v_i^2 \rangle =-2,\\
\langle v_i,v_j \rangle \geq 0$, $i \ne j,\\
\langle v,v_i \rangle=\langle \widehat{H},v_i \rangle=0.
\end{cases}
\end{equation}
If there is a decomposition
$\{0,1,\dots,n\}=I \amalg J$
such that $\langle v_i,v_j \rangle =0$ for all $i \in I, j \in J$,
then $0=\langle v^2 \rangle=\langle(\sum_{i \in I}a_i v_i)^2 \rangle
+\langle (\sum_{j \in J}a_j v_j)^2 \rangle$.
By \eqref{eq:hodge2}, $\sum_{i \in I}a_i v_i=\sum_{j \in J}a_j v_j=0$
in $v^{\perp}/{\Bbb Z}v$.
Since $0 \leq \rk (\sum_{i \in I}a_i v_i), 
\rk (\sum_{j \in J}a_j v_j)\leq \rk v$,
we get $I=\emptyset$ or $J=\emptyset$.
Then as in the classification of the singular fiber of an
elliptic surface, we get our claim.
\end{proof}


Assume that there is another point $x'$ of $\overline{M}_H^{v}(v)$ 
representing a properly $v$-twisted semi-stable sheaf
$\bigoplus_{i=0}^{n'} {E_i'}^{\oplus a_i'}$, where $E_i'$
is a $v$-twisted stable sheaf such that
$\langle v(E_i'),\widehat{H} \rangle=\langle v(E_i'),v \rangle=0$
and $E_i' \ne E_j'$ for $i \ne j$.  We set $v_i':=v(E_i')$. 
Then we have the following lemma. 
\begin{lem}\label{lem:perp}
\begin{equation}
({\Bbb Z}v_0+{\Bbb Z}v_1+\dots+{\Bbb Z}v_n) \perp 
({\Bbb Z} v_0'+{\Bbb Z}v_1'+\dots+{\Bbb Z}v_{n'}').
\end{equation}
\end{lem}

\begin{proof}
We set
\begin{equation}
\begin{split}
S_1:=&\{\;i\; |\; 
v_i' \in ({\Bbb Z}v_0+{\Bbb Z}v_1+\dots+{\Bbb Z}v_n)^{\perp} \},\\
S_2:=& \{0,1,\dots,n'\}\setminus S_1.
\end{split}
\end{equation}
Then $v=\sum_{i \in S_1}a_i' v_i'+\sum_{i \in S_2}a_i' v_i'$.
Assume that $i \in S_2$. Since $0=\langle v_i',v \rangle=
\sum_j a_j\langle v_i',v_j \rangle$,
$\langle v_i',v_j \rangle<0$ for some $j$.
Then $\chi(E_i',E_j)>0$, which implies that
there is a non-zero homomorphism
$E_i' \to E_j$ or $E_j \to E_i'$.
Since $E_i'$ and $E_j$ are $v$-twisted stable sheaves such that
$\chi(E,E_i'(nH))/\rk E_i'=\chi(E,E_j(nH))/\rk E_j
(=\chi(E,E(nH))/\rk E)$ for all $n$, we get
$E_i' \cong E_j$.
Thus $v_i' \in \{v_0,v_1,\dots,v_n \}$. 
Then we get that
$\langle (\sum_{i \in S_1}a_i' v_i')^2 \rangle+
\langle(\sum_{i \in S_2}a_i' v_i')^2 \rangle=0$,
%
%
%
%
and hence $\sum_{i \in S_1}a_i' v_i', \sum_{i \in S_2}a_i' v_i' \in 
{\Bbb Z}v$.
Since $\rk(\sum_{i \in S_2}a_i' v_i') \leq \rk v$,
$\sum_{i \in S_2}a_i' v_i'=0$ or $v$.
If $\sum_{i \in S_2}a_i' v_i'=v$, then 
Lemma \ref{lem:lattice} implies that
$S_2=\{0,1,\dots, n \}$ and $a_i'=a_i$, which implies that $x'=x$.
Since $x \ne x'$, we get that $S_1=\{0,1,\dots,n'\}$.
Thus our claim holds.
\end{proof}

We shall study the fiber of $\phi_{\alpha}:M_H^{v+\alpha}(v) \to
\overline{M}_H^v(v)$ at $x$.
By the classification of the (extended) Dynkin diagram, 
we may assume that $a_0=1$.
Then $v_i$, $i=1,2,\dots,n$ become a fundamental root system
of the corresponding finite Lie algebra $\frak g$
under the change of the sign of the bilinear form.

\begin{lem}
\begin{equation}
{\cal U}'=({\cal U}' \cap (\bigoplus_{i=0}^n {\Bbb Z} v_i)^{\perp})
\coprod ({\cal U}'\cap \bigoplus_{i=0}^n {\Bbb Z} v_i)
\end{equation}
and
\begin{equation}
\begin{split}
{\cal U}'\cap (\bigoplus_{i=0}^n {\Bbb Z} v_i)
&=\{ u \in \bigoplus_{i=0}^n {\Bbb Z} v_i| \langle u^2 \rangle=-2, 0<\rk u <\rk v \}\\
&=\Psi_+ \coprod (v-\Psi_+),
\end{split}
\end{equation}
where $\Psi_+:=\{ u=\sum_{i=1}^n b_i v_i| \langle u^2 \rangle=-2,
b_i \geq 0 \}$ is the set of positive roots of ${\frak g}$.

\end{lem}

\begin{proof}
For $u \in {\cal U}'$, we set $w=v-u$.
Since $\langle v,u \rangle=0$,
we get $\langle w^2 \rangle=\langle u^2 \rangle=-2$.
Since $0<\rk u<\rk v$, we have $\rk w>0$.
By Proposition \ref{prop:exist},
there are $v$-twisted semi-stable sheaves
$F$ and $G$ with $v(F)=u$ and $v(G)=w$.
Applying Lemma \ref{lem:perp} to $F \oplus G$, we see that
$u=\sum_{i=0}^n b_i v_i \in \bigoplus_{i=0}^n {\Bbb Z} v_i,
a_i \geq b_i \geq 0$ or
$u \in (\bigoplus_{i=0}^n {\Bbb Z} v_i)^{\perp}$
according as $F \oplus G$ is $S$-equivalent to $\bigoplus_{i=0}^n E_i^{\oplus a_i}$
 or not.
Thus the first claim holds.
If $b_0=0$, then $u \in \Psi_+$ and if $b_0=1$, then $w \in \Psi_+$.
Thus $u=v-w \in v-\Psi_+$. Therefore the second assertion also holds.
\end{proof}
%
%
Therefore the wall $W_u$ corresponds to the wall
defining the Weyl chamber.
More precisely, $W_u \cap ((\bigoplus_{i=0}^n {\Bbb Z} v_i)/{\Bbb Z}v)
\otimes {\Bbb R}$ is the corresponding wall.
We define the fundamental Weyl chamber:

\begin{equation}
D:=\{\alpha \in \delta(H^{\perp})_{alg} \otimes_{\Bbb Q} {\Bbb R}\;|
\;\langle v_i,\alpha \rangle>0, i>0\}.
\end{equation}
For a small $\alpha \in D$, we describe the exceptional set
$\phi_{\alpha}^{-1}(x)$.
The method is the same as in \cite{Y:5}.

\begin{lem}\label{lem:basic1}
Assume that $\alpha \in \delta(H^{\perp})_{alg}$ belongs to $D$
and $|\langle \alpha^2 \rangle| \ll 1$.
Let $F$ be a $v+\alpha$-twisted semi-stable sheaf such that
$v(F)=v_0+\sum_{j>0}b_j v_j$,
$0 \leq b_j \leq a_j$.
\begin{enumerate}
\item[(1)]
If $v(F) \ne v$, then
$F$ is $v+\alpha$-twisted stable and $F$ is $S$-equivalent to 
$E_0 \oplus (\bigoplus_{j>0} E_j^{\oplus b_j})$ with respect to
the $v$-twisted stability.
\item[(2)]
For a non-zero homomorphsim $\phi:E_i \to F$, $i>0$,
$\phi$ is injective and
$F':=\coker \phi$ is a $v+\alpha$-twisted stable sheaf.
\item[(3)]
If there is a non-trivial extension
\begin{equation}\label{eq:ext4}
0 \to E_i \to F'' \to F \to 0
\end{equation}
and $b_i+1 \leq a_i$,
then $F''$ is $v+\alpha$-twisted stable.
\end{enumerate}
\end{lem}

\begin{proof}
We take elements $E,A \in K(X) \otimes {\Bbb Q}$ such that
$v(E)=v$, $v(A)=\alpha$.
Since $|\langle \alpha^2 \rangle| \ll 1$, $F$ is $v$-twisted semi-stable.
Assume that $F$ is $S$-equivalent to 
$\bigoplus_{j=0}^{n'} {E_j'}^{\oplus b_j'}$ with respect to
the $v$-twisted stability.
Since $F \oplus (\bigoplus_{j>0} 
E_j^{\oplus (a_j-b_j)})$ is $S$-equivalent to
$\bigoplus_{j=0}^{n'} {E_j'}^{\oplus b_j'} \oplus (\bigoplus_{j>0} 
E_j^{\oplus (a_j-b_j)})$, 
by Lemma \ref{lem:lattice} and
Lemma \ref{lem:perp}, we get that $F$ is $S$-equivalent to
$E_0 \oplus (\bigoplus_{j>0} E_j^{\oplus b_j})$ with respect to
the $v$-twisted stability.
Since $\chi(E+A,F) \geq \chi(E+A,E)=0$ and $\chi(E+A,E_i)<0$ for all $i>0$,
there is no proper subsheaf $E'$ such that
$\chi(E+A,E'(nH))/\rk E'=\chi(E+A,F(nH))/\rk F$ for all $n$. 
Thus the claim (1) holds.  

We next prove (2).
Since $E_i$ is $v$-twisted stable and
$F$ is $v$-twisted semi-stable,
$\phi$ is injective and
$F'$ is a $v$-twisted semi-stable sheaf.
By (1), $F'$ is $S$-equivalent to
$E_0 \oplus (\bigoplus_{j>0} E_j^{\oplus c_j})$
with respect to the $v$-twisted stability, where
$v(F')=v_0+\sum_{j>0}c_j v_j$. 
If $F'$ is not $v+\alpha$-twisted stable, then 
there is a quotient sheaf $F' \to G$ such that
$\deg_E G=0$ and 
$\chi(E+A,G)/\rk G<\chi(E+A,F')/\rk F'$.
Since $|\langle \alpha^2 \rangle| \ll 1$ and $\chi(E+A,E_i)<0$
for all $i>0$, we see that
$G$ is $v$-twisted semi-stable and is 
$S$-equivalent to $\bigoplus_{j>0} E_j^{\oplus c_j'}$
with respect to the $v$-twisted stability.
Hence we get that $\chi(E+A,G)<0$, which implies that $F$ is not
$v+\alpha$-twisted stable.
Therefore $F'$ is $v+\alpha$-twisted stable.

Finally we prove (3).
By our assumption,
$\chi(E+A,F'') \geq \chi(E+A,E)=0$.
If $F''$ is not $v+\alpha$-twisted stable, then
there is a quotient sheaf $F'' \to G$ such that
$\deg_E G=0$ and $\chi(E+A,G)/\rk G<\chi(E+A,F'')/\rk F''$.
Then we see that
$G$ is $v$-twisted semi-stable and is 
$S$-equivalent to $\bigoplus_{j>0} E_j^{\oplus c_j''}$
with respect to the $v$-twisted stability.
Then there is a quotient $G \to E_j$, $j>0$.
By \eqref{eq:ext},
we have an exact sequence
\begin{equation}
0=\Hom(F,E_j) \to \Hom(F'',E_j) \to \Hom(E_i,E_j).
\end{equation}
We consider the map $\psi:E_i \to F'' \to G \to E_j$.
Then $\psi$ is an isomorphism, which implies that 
the extension \eqref{eq:ext4} splits.  
\end{proof}

By Lemma \ref{lem:basic1} (1), we get the following corollary.

\begin{cor}\label{cor:stable}
We set $w:=v_0+\sum_{j>0}b_j v_j$, $0 \leq b_j \leq a_j$.
If $w \ne v$, then $\overline{M}_H^{v+\alpha}(w)=M_H^{v+\alpha}(w)$.
In particular, if $\langle w^2 \rangle=-2$, then $M_H^{v+\alpha}(w)$
is not empty and consists of one element.
\end{cor}

\begin{cor}\label{cor:dim}
Let $F$ be a $v+\alpha$-twisted stable sheaf with 
$v(F)=v_0+\sum_{j>0}b_j v_j$, $0 \leq b_j \leq a_j$.
\begin{enumerate}
\item[(1)]
If $v(F)=v$, then
$\dim \Hom(E_i,F) \leq 1$.
\item[(2)]
If $v(F) \ne v$, then 
$\dim \Hom(E_i,F)=\max\{-\langle v(F), v_i \rangle,0 \}$.
\end{enumerate}
\end{cor}

\begin{proof}
We set $\dim \Hom(E_i,F)=k$.
By the Riemann-Roch theorem,
$k \geq -\langle v(F),v_i \rangle$.
Hence if $k=0$, then 
our claims (1), (2) hold.
Assume that $k>0$.
By Lemma \ref{lem:basic1},
$\phi:\Hom(E_i,F) \otimes E_i \to F$ is injective
and $F':=\coker \phi$ is $v+\alpha$-twisted stable.
If $v(F)=v$, then
$-2 \leq \langle v(F')^2 \rangle =-2k^2$.
Hence $k \leq 1$.
If $v(F) \ne v$, then $\langle v(F)^2 \rangle=-2$ and hence
$-2 \leq \langle v(F')^2 \rangle =-2-2k(k+\langle v(F),v_i \rangle)$.
Then $k+\langle v(F),v_i \rangle \leq 0$, which implies that
$k=- \langle v(F),v_i \rangle$. 
\end{proof}

\begin{cor}\label{cor:par1}
We set $w=v_0+\sum_{j>0} b_j v_j$, $0 \leq b_j \leq a_j$.
If $\langle w^2 \rangle=\langle (w-v_i)^2 \rangle=-2$,
then we have an isomorphism 
$M_H^{v+\alpha}(w) \to M_H^{v+\alpha}(w-v_i)$
sending $F$ to $\coker(E_i \to F)$.
\end{cor}

\begin{proof}
By our assumption, we see that $\langle w,v_i \rangle=-1$.
By Corollary \ref{cor:dim},
$\Hom(E_i,F)={\Bbb C}$.
By Lemma \ref{lem:basic1} (2), 
$F':=\coker(E_i \to F)$ is a $v+\alpha$-twisted stable sheaf
with $v(F')=w-v_i$. 
Conversely for a $v+\alpha$-twisted stable sheaf
$F'$ with $v(F')=w-v_i$, we get $\langle w-v_i,v_i \rangle=1$,
and hence by Corollary \ref{cor:dim} and 
Lemma \ref{lem:basic1} (3), the non-trivial extension of
$F'$ by $E_i$ gives a $v+\alpha$-twisted stable sheaf $F$ with
$v(F)=w$.
\end{proof}

We set
\begin{equation}
\begin{split}
C_i:=&\left\{(E,U)\left|
\begin{split}
&E \in M_H^{v+\alpha}(v), U \subset \Hom(E_i,E)\\
& \dim U=1
\end{split} \right. \right\}\\
=&\{E_i \subset E| E \in M_H^{v+\alpha}(v) \}.
\end{split}
\end{equation}
$C_i$ is the moduli space of twisted coherent systems.
\begin{prop}\label{prop:BN}
\begin{enumerate}
\item[(1)]
$C_i \subset \phi_{\alpha}^{-1}(x)$.
\item[(2)]
$C_i \cong {\Bbb P}^1$ and the natural map
$\pi:C_i \to M_H^{v+\alpha}(v)$ is a closed immersion.
In particular, $C_i$ is not empty.
\end{enumerate}
\end{prop}

\begin{proof}
We set $F:=\coker(E_i \to E)$.
Then $E$ is $S$-equivalent to $E_i \oplus F$ with respect to the
$v$-twisted stability.
By Lemma \ref{lem:basic1} (1),
$F$ is $S$-equivalent to $E_i^{\oplus (a_i-1)}\oplus 
\bigoplus_{j \ne i} E_j^{\oplus a_j}$, and hence the first claim holds.  
We next show the assertion (2).
We note that the Zariski tangent space of $C_i$ at
$E_i \to E$ is
\begin{equation}\label{eq:tangent}
 \Ext^1(E_i \to E,E)=\Ext^1(F,E)
\end{equation}
and the obstruction for the infinitesimal lifting belongs to
\begin{equation}\label{eq:obstruction}
\ker(\Ext^2(E_i \to E,E) \to \Ext^2(E,E) \overset{tr}{\to} 
H^2(X,{\cal O}_X)).
\end{equation}
We shall first show that $C_i$ is smooth at $E_i \subset E$.
Since $\Ext^2(E_i \to E,E)=\Ext^2(F,E)=\Hom(E,F)^{\vee}$,
it is sufficient to show that
$\Hom(E,F)={\Bbb C}$.
By the exact sequence
\begin{equation}
 0 \to E_i \to E \to F \to 0,
\end{equation}
we get an exact sequence
\begin{equation}
 0 \to \Hom(F,F) \to \Hom(E,F) \to \Hom(E_i,F).
\end{equation}
If $\Hom(E_i,F) \ne 0$, then we get
$\dim \Hom(E_i,E) \geq 2$, which contradicts Corollary \ref{cor:dim}.
Hence $ \Hom(F,F) \cong \Hom(E,F)$.
By Lemma \ref{lem:basic1} (2),
$F$ is simple. Therefore $\Hom(E,F)={\Bbb C}$.
Since the homomorphism
\begin{equation}
 \Ext^1(E_i \to E,E) \to \Ext^1(E,E)
\end{equation}
between the Zariski tangent spaces is injective,
$C_i \to M_H^{v+\alpha}(v)$ is a closed immersion, 
provided that $C_i \ne \emptyset$.

We next show that $C_i \ne \emptyset$ and isomorphic to
${\Bbb P}^1$.
Since $\langle (v-v_i)^2 \rangle =-2$,
Corollary \ref{cor:stable} implies that
$M_H^{v+\alpha}(v)$ consists of exactly one 
$v+\alpha$-twisted stable sheaf $F$.
By Corollary \ref{cor:dim} (2), $\Hom(E_i,F) =0$. 
Thus $\Ext^2(F,E_i)=\Hom(E_i,F)^{\vee}=0$. 
Since $\Hom(F,E_i)=0$, we get 
$\Ext^1(F,E_i) \cong {\Bbb C}^{\oplus 2}$. 
Let $E$ be a coherent sheaf which is defined by a non-trivial extension
\begin{equation}\label{eq:ext}
 0 \to E_i \to E \to F \to 0.
\end{equation}
By Lemma \ref{lem:basic1} (3), $E$ is $v+\alpha$-twisted stable.
Therefore $C_i \ne \emptyset$ and $C_i \cong {\Bbb P}^1$. 
\end{proof}

\begin{prop}\label{prop:config}
We identify $C_i$ with its image $\pi(C_i)$.
Then we have 
\begin{enumerate}
\item[(1)]
\begin{equation}
 \phi_{\alpha}^{-1}(x)=\cup_{i=1}^n C_i
\end{equation}
and $\cup_{i=1}^n C_i$ is a simple normal crossing divisor.
\item[(2)]
$(C_i,C_j)=\langle v_i,v_j \rangle$.
In particular the dual graph of $C_i$, $1 \leq i \leq n$
is of type $A_n,D_n,E_n$. 
\end{enumerate}
\end{prop}

\begin{proof}
By Proposition \ref{prop:BN} (1),
$\phi_{\alpha}^{-1}(x)\supset \cup_{i=1}^n C_i$.
Since $\alpha \in D$, we get $\phi_{\alpha}^{-1}(x)=\cup_{i=1}^n C_i$.
We shall study the configuration of $C_i$, $1 \leq i \leq n$.
Here we give a geometric argument based on Lemma \ref{lem:basic1}. 
We shall prove the following assertions:
\begin{enumerate}
\item
$C_i \cap C_j \ne \emptyset$ if and only if $\langle v_i,v_j \rangle=1$.
\item
If $C_i \cap C_j \ne \emptyset$, then $\#(C_i \cap C_j)=1$.
\item
If $C_i \cap C_j \ne \emptyset$, then $C_i$ and $C_j$ intersect
transversely.
\item
$C_i \cap C_j \cap C_k= \emptyset$ for three curves $C_i,C_j,C_k$. 
\end{enumerate}

(i) Assume that $C_i \cap C_j \ne \emptyset$
and take a point $E \in C_i \cap C_j$.
Then $E$ fits in an exact sequence
\begin{equation}\label{eq:ext2}
 0 \to E_i \oplus E_j \to E \to F' \to 0
\end{equation}
By Lemma \ref{lem:basic1} (2),
$F'$ is a $v+\alpha$-twisted stable sheaf.
Since $\delta(H^{\perp})_{alg}$ is negative definite,
$-2 \leq \langle v(F')^2 \rangle <0$. Thus $\langle v(F')^2 \rangle=-2$. Then
$-2 =\langle v(F')^2 \rangle=
-4+2\langle v(E_i),v(E_j) \rangle$, which implies that
$\langle v(E_i),v(E_j) \rangle=1$.
Conversely if $\langle v(E_i),v(E_j) \rangle=1$, then
$\langle (v-(v_i+v_j))^2 \rangle=-2$. 
Hence there is a $v+\alpha$-twisted semi-stable sheaf $F'$
with $v(F')=v-(v_i+v_j)$.
By Lemma \ref{lem:basic1} (1), $F'$ is $v+\alpha$-twisted stable.
By Corollary \ref{cor:dim} (2), $\Hom(E_i,F')=0$. 
Hence $\dim \Ext^1(E_i,F')=\langle v_i,v-(v_i+v_j) \rangle=1$.
We also have $\dim \Ext^1(E_j,F')=1$.
We take an extension 
\begin{equation}\label{eq:ext3}
 0 \to E_i \oplus E_j \to E \to F' \to 0
\end{equation}
whose extension class is given by $(e_i,e_j) \in
\Ext^1(E_i,F') \oplus \Ext^1(E_j,F')$, $e_i, e_j \ne 0$.
Then Lemma \ref{lem:basic1} (3) implies that $E$ is 
a $v+\alpha$-twisted stable sheaf with
$v(E)=v$. 
Therefore $C_i$ and $C_j$ intersect at $E$.

(ii) Assume that $C_i \cap C_j \ne \emptyset$.
Then every member of $C_i \cap C_j$ 
fits in an extension \eqref{eq:ext3}.
Since $M_H^{v+\alpha}(v-(v_i+v_j))=\{F'\}$ 
and $E$ does not depend on the choice of $(e_i,e_j)$,
we get $\#(C_i \cap C_j)=1$.

(iii)
Assume that $C_i$ and $C_j$ intersect at $E$. 
Since $\Hom(E_i \oplus E_j,F')=0$ and $F'$ is $v+\alpha$-twisted stable,
$\Ext^2(F',E)=\Hom(E,F')^{\vee}={\Bbb C}$.
Then we see that the natural homomorphism
\begin{equation}
 \Ext^1(E/E_i,E) \oplus \Ext^1(E/E_j,E) \to \Ext^1(E,E)
\end{equation}
of tangent spaces 
is an isomorphism, and hence 
$C_i$ and $C_j$ intersect transversely.

(iv)
If $C_i \cap C_j \cap C_k \ne \emptyset$, then
$\langle v_i,v_j \rangle=\langle v_j,v_k \rangle=\langle v_k,v_i \rangle=1$,
which implies that $\langle (v_i+v_j+v_k)^2 \rangle=0$.
Since $\delta(H^{\perp})_{alg}$ is negative definite,
this is impossible.

Therefore $\cup_i C_i$ is simple normal crossing and
$(C_i,C_j)=\langle v_i,v_j \rangle$.
\end{proof}

\begin{lem}\label{lem:equation}
Let $[C_i] \in H_2(M_H^{v+\alpha}(v), {\Bbb Z})$ 
be the fundamental class of $C_i$ and
$\PD([C_i]) \in H^2(M_H^{v+\alpha}(v), {\Bbb Z})$ the 
Poincar\'{e} dual of $[C_i]$.
Then
\begin{equation}\label{eq:C_i}
 \PD([C_i])=\theta_v^{\alpha}(-v_i).
\end{equation}
\end{lem}

\begin{proof}
Let $M_H^{v+\alpha}(v) =\cup_{\lambda} U_{\lambda}$ be 
an analytic open covering and ${\cal F}_{\lambda}$ 
a local universal family on
$U_{\lambda} \times X$.
Then $\Hom_{p_{U_{\lambda}}}(p_X^*(E_i),{\cal F}_{\lambda})=
\Ext^2_{p_{U_{\lambda}}}(p_X^*(E_i),{\cal F}_{\lambda})=0$ and 
$C_i \cap U_{\lambda}$ is the scheme-theoretic support of
$\Ext^1_{p_{U_{\lambda}}}(p_X^*(E_i),{\cal F}_{\lambda})$, 
where $p_X:U_{\lambda} \times X \to X$ and
$p_{U_{\lambda}}:U_{\lambda} \times X \to U_{\lambda}$ are projections.
For a sufficiently large integer $n$,
$V_{\lambda}:=
p_{U_{\lambda}*}({\cal F}_{\lambda} \otimes p_X^*({\cal O}_X(nH)))$
is a locally free sheaf on $U_{\lambda}$.
Then we can glue
$\{{\cal F}_{\lambda} \otimes p_{U_{\lambda}}^*(V_{\lambda}^{\vee})\}$
together and we get a quasi-universal family ${\cal F}$ on
$M_H^{v+\alpha}(v) \times X$.
By using the Grothendieck Riemann-Roch theorem, we get that
the Poincar\'{e} dual of $[C_i]$ is $\theta_v^{\alpha}(-v_i)$.
\end{proof}

\begin{rem}
By Lemma \ref{lem:equation},
the non-emptyness of $C_i$ also follows from the fact that 
$\theta_v^{\alpha}(v_i) \ne 0$.
Since $\theta_v^{\alpha}$ is an isometry,
the configuration of $C_i$ is also 
described by the configuration of $v_i$.
In particular, we get a different proof of Proposition \ref{prop:config}. 
\end{rem} 

Since $\phi_{\alpha}^{-1}(x)$ is a union of $(-2)$-curves, 
we get the following proposition.

\begin{prop}
$\phi_{\alpha}:M_H^{v+\alpha}(v) \to \overline{M}_H^v(v)$
is surjective and 
$\overline{M}_H^v(v)$ contains a $v$-twisted stable sheaf.
\end{prop}

\begin{proof}
Let $x$ be a point of $\overline{M}_H^v(v)$ and assume that
$x$ corresponds to a properly $v$-twisted semi-stable sheaf.
Then $\phi_{\alpha}^{-1}(x) \ne \emptyset$ with
$\dim \phi_{\alpha}^{-1}(x)=1$.
Since $\phi_{\alpha}:\phi_{\alpha}^{-1}(M_H^v(v)) \to M_H^v(v)$
is an isomorphism and $\overline{M}_H^v(v) \setminus M_H^v(v)$
is a finite set, we get our claims.
\end{proof}

%

The remaining of this section is an appendix.
\begin{lem}\label{lem:basic2}
Assume that $\alpha \in D$, $|\langle \alpha^2 \rangle| \ll 1$
 satisfies that
\begin{equation}\label{eq:slope}
 \langle v_i,\alpha \rangle/\rk v_i> 
 \langle v+\sum_{j>0} a_j v_j,\alpha \rangle/\rk(v+\sum_{j>0} a_j v_j)
\end{equation}
for all $i>0$.
Let $F$ be a $v+\alpha$-twisted semi-stable sheaf such that
$v(F)=v+\sum_{j>0} b_j v_j$,
$0 \leq b_j \leq a_j$. Then
\begin{enumerate}
\item[(1)]
$F$ is $v+\alpha$-twisted stable.
\item[(2)]
Let $F'$ be a coherent sheaf which fits in a non-trivial extension
\begin{equation}
 0 \to E_i \to F' \to F \to 0
\end{equation} 
and $b_i+1 \leq a_i$.
Then $F'$ is $v+\alpha$-twisted stable.
\item[(3)]
For a subsheaf $E_i$, $i>0$ of $F$,
$F':=F/E_i$ is $v+\alpha$-twisted stable.
\end{enumerate}

\end{lem}

\begin{proof}
We take $E$, $A \in K(X) \otimes {\Bbb Q}$ with
$v(E)=v$, $v(A)=\alpha$.
Let 
\begin{equation}
 0=F_0 \subset F_1 \subset F_2 \subset \dots \subset F_s=F
\end{equation}
be the Jordan-H\"{o}lder filtration of $F$ with respect to the
$v$-twisted stability.
For $F_i/F_{i-1}$ with $F_i/F_{i-1} \not \in
\{E_0,E_1,\dots, E_n \}$, a similar argument to the proof of
Lemma \ref{lem:perp} shows that
$\Ext^1(F_i/F_{i-1},E_j)=0$ for all $j$.
Hence replacing the filtration, 
we may assume that $F_i/F_{i-1} \not \in  \{E_0,E_1,\dots, E_n \}$
for $0 \leq i \leq k$ and
$F_i/F_{i-1} \in \{E_0,E_1,\dots, E_n \}$ for $i>k$.
Then we get $F=F_k \oplus F/F_k$ and
$-2 \leq \langle v(F)^2 \rangle =
\langle v(F_k)^2 \rangle +\langle v(F/F_k)^2 \rangle$.
Assume that $F_k \ne 0$, that is, $k>0$.
Then we see that (i) $v(F_k)=v$ and 
$\langle v(F/F_k)^2 \rangle=-2$, or (ii) 
$\langle v(F_k)^2 \rangle=-2$ and 
$v(F/F_k) \in {\Bbb Z}v$.
If $v(F_k)=v$, then $\langle v(F_k),\alpha \rangle=0$ and
$\langle v(F/F_k),\alpha \rangle=\sum_j b_j \langle v_j,\alpha \rangle >0$, 
which contradicts the 
$v+\alpha$-twisted semi-stability of $F$.
Hence the case (i) does not occur.
If the case (ii) occurs, then 
since $\langle v(F_k)^2 \rangle=-2$,
we get that
$\langle v(F),v(F_i/F_{i-1}) \rangle=
\langle v(F_k),v(F_i/F_{i-1}) \rangle \ne 0$ for some $i \leq k$.
On the other hand,
by our choice of $v(F)$, we get 
$\langle v(F), v(F_i/F_{i-1}) \rangle=0$ for all $i \leq k$.
Hence the case (ii) does not occur.
Therefore $F_k=0$. 
Then we see that $F$ is $S$-equivalent to
$E_0 \oplus (\bigoplus_{j>0} E_j^{\oplus (a_j+b_j)})$.

Our assumption \eqref{eq:slope} implies that 
\begin{equation}
 \chi(E+A,E_i)/\rk E_i<\chi(E+A,2E-E_0)/\rk(2E-E_0).
\end{equation}
Since $F=2E-E_0-\sum_{j>0} c_j E_j$, $0 \leq c_j \leq a_j$
as an element of $K(X)$,
we get that 
\begin{equation}
 \chi(E+A,2E-E_0)/\rk (2E-E_0) \leq 
 \chi(E+A,F)/\rk F.
\end{equation}
In the same way as in the proof of Lemma \ref{lem:basic1} (1),
we see that
$F$ is a $v+\alpha$-twisted stable sheaf.
Thus (1) holds.
The proof of (2) and (3) are the same as in the proof of 
Lemma  \ref{lem:basic1}.
\end{proof}

\begin{rem}
If $\langle v_i,\alpha \rangle=\langle v_j,\alpha \rangle$ for all
$1 \leq i,j \leq n$,
then \eqref{eq:slope} is satisfied.
\end{rem}   

In the same way as in the proof of Corollary \ref{cor:par1},
we get the following.
\begin{cor}\label{cor:par2}
Assume that $\alpha \in D$, $|\langle \alpha^2 \rangle| \ll 1$
 satisfies \eqref{eq:slope}.
We set $w=v+\sum_{j>0} b_j v_j$, $0 \leq b_j \leq a_j$.
\begin{enumerate}
\item[(1)]
If $w \ne v$, then $\overline{M}_H^{v+\alpha}(w)=M_H^{v+\alpha}(w)$.
\item[(2)]
If $\langle w^2 \rangle=\langle (w-v_i)^2 \rangle=-2$,
then we have an isomorphism 
$M_H^{v+\alpha}(w) \to M_H^{v+\alpha}(w-v_i)$
sending $F \in M_H^{v+\alpha}(w)$ to $\coker(E_i \to F)$.
\end{enumerate}
\end{cor}

\begin{rem}
Assume that $\alpha \in D$, $|\langle \alpha^2 \rangle| \ll 1$
 satisfies \eqref{eq:slope}.
We note that $H_2(\phi_{\alpha}^{-1}(x),{\Bbb C}) \to
H_2(M_H^{v+\alpha}(v),{\Bbb C})$ is injective.
We can regard $H_2(\phi_{\alpha}^{-1}(x),{\Bbb C})$ as the Cartan subalgebra
of ${\frak g}$.
In order to get ${\frak g}$,
we set 
\begin{equation}
 \Psi:=\{u| u=\sum_{i=1}^n b_i v_i, 
 \langle u^2 \rangle=-2\}.
\end{equation}
Let $P(w,w-v_i)$
be the subscheme of $M_H^{v+\alpha}(w) \times M_H^{v+\alpha}(w-v_i)$
consisting of pairs $(E,F) \in M_H^{v+\alpha}(w) \times M_H^{v+\alpha}(w-v_i)$
which fits in an exact sequence
\begin{equation}
0 \to E_i \to E \to F \to 0.
\end{equation}
Then we can show that $P(w,w-v_i)$ is isomorphic to ${\Bbb P}^1$ or a point.
As in \cite{Na:1}, we see that there is an action of ${\frak g}$ on
$H_2(M_H^{v+\alpha}(v),{\Bbb C}) \oplus \bigoplus_{u \in \Psi} 
H_0(M_H^{v+\alpha}(v+u),{\Bbb C})$
and we have an isomorphism of 
${\frak g}$-module:
\begin{equation}
H_2(\phi_{\alpha}^{-1}(x),{\Bbb C}) \oplus \bigoplus_{u \in \Psi} 
H_0(M_H^{v+\alpha}(v+u),{\Bbb C}) \cong {\frak g}.
\end{equation}
For a homology class $[x] \in 
H_{2+\langle w^2 \rangle}(M_H^{v+\alpha}(w),{\Bbb C})$,
the action of Chevalley generators $e_i,f_i,h_i$, $1 \leq i \leq n$ 
are given by
\begin{equation}
\begin{split}
e_i:[x] \mapsto & p_{M_H^{v+\alpha}(w+v_i)*}(
(M_H^{v+\alpha}(w+v_i) \times [x]) \cap P(w+v_i,w))\\
f_i:[x] \mapsto & (-1)^{t(w)} p_{M_H^{v+\alpha}(w-v_i)*}(
([x] \times M_H^{v+\alpha}(w-v_i)) \cap P(w,w-v_i))\\
h_i:[x] \mapsto & -\langle w,v_i \rangle [x],
\end{split}
\end{equation}
where $p_{M_H^{v+\alpha}(w+kv_i)}:M_H^{v+\alpha}(w) \times
M_H^{v+\alpha}(w-v_i) \to M_H^{v+\alpha}(w+kv_i)$, $k=0,-1$ are projections
and $t(w)=(\dim M_H^{v+\alpha}(w-v_i)-\dim M_H^{v+\alpha}(w))/2
=-(\langle w, v_i \rangle+1)$.

%
\end{rem}

\subsection{Other chambers}

\begin{defn}\label{defn:weyl}
Let $W$ be the Weyl group generated by reflections
$R_{v_i}$, $i=1,2,\dots,n$.
\end{defn}
$W$ is the Weyl group of ${\frak g}$.
By Lemma \ref{lem:equation} and
Proposition \ref{prop:comm} (also see Remark \ref{rem:comm}), 
we get the following. 
\begin{prop}
If $\alpha \in w(D)$, $w \in W$, then
$\PD([C_i])=\theta_v^{\alpha}(-w(v_i))$.
\end{prop}
We shall give a set-theoretic description of $C_i$.
Let $M_H^{v+\alpha}(v) =\cup_{\lambda} U_{\lambda}$ be 
an analytic open covering and ${\cal F}_{\lambda}$ 
a local universal family on
$U_{\lambda} \times X$.
If $\rk (w(v_i))=0$, then we set $w(v_i):=\eta+b \rho_X$,
where $\eta \in \Pic(X)$ satisfies $(\eta^2)=-2$.
Then $0= \langle \widehat{H},w(v_i) \rangle=(H,\eta)$, which 
contradicts the ampleness of $H$.
Hence $\rk w(v_i) \ne 0$.
Let $F_i$ be a $v+\alpha$-twisted semi-stable sheaf with
$v(F_i)=\pm w(v_i)$.
We first assume that $\rk (w(v_i))>0$, that is,
$v(F_i)=w(v_i)$.
Since $\langle \alpha, w(v_i) \rangle=
\langle w^{-1}(\alpha), v_i \rangle > 0$,
$\chi(E+A,F_i)=-\langle v+\alpha,v(F_i) \rangle=
-\langle \alpha,w(v_i) \rangle<0$.  
Then $\Ext^2(F_i,E)=0$ for all $E \in M_H^{v+\alpha}(v)$.
If $E$ is $v$-twisted stable, then $\Hom(F_i,E)=0$.
Hence $L_{\lambda}:=\Ext^1_{p_{U_{\lambda}}}(p_X^*(F_i),{\cal F}_{\lambda})$ 
is a torsion sheaf of pure dimension 1
whose support is contained in $\phi_{\alpha}^{-1}(x)$.
By the Grothendieck Riemann-Roch theorem, the Poincar\'{e} dual 
of the scheme-theoretic support of $L_{\lambda}$ is 
$\theta_v^{\alpha}(-w(v_i))=\PD([C_i])$.
Since $H^0(M_H^{v+\alpha}(v),{\cal O}_{M_H^{v+\alpha}(v)}(C_i))={\Bbb C}$,
 we get that
\begin{equation}
\begin{split}
C_i&=\{E \in M_H^{v+\alpha}(v)|\Ext^1(F_i,E) \ne 0 \}\\
&=\{E\in M_H^{v+\alpha}(v)|\Hom(F_i,E) \ne 0 \}.
\end{split}
\end{equation}
If $\chi(E+A,F_i)>0$, that is, $\rk (w(v_i))<0$, then
we see that $\Hom(F_i,E)=0$ for all $E \in M_H^{v+\alpha}(v)$
and $\Ext^2(F_i,E)=0$ for a $v$-twisted stable sheaf $E$.
Then we also have $\Ext^1(F_i,E)=0$. Hence we see that
$\Ext^1_{p_{U_{\lambda}}}(p_X^*(F),{\cal F}_{\lambda})=0$ and 
$\Ext^2_{p_{U_{\lambda}}}(p_X^*(F),{\cal F}_{\lambda})$ 
is a torsion sheaf of pure dimension 1. 
Hence we get that
\begin{equation}
\begin{split}
C_i&=\{E \in M_H^{v+\alpha}(v)|\Ext^2(F_i,E) \ne 0 \}\\
&=\{E\in M_H^{v+\alpha}(v)|\Hom(E,F_i) \ne 0 \}.
\end{split}
\end{equation}
Therefore we get the following proposition.
\begin{prop}
Assume that $\alpha \in w(D)$, $w \in W$.
Then $w(v_i) \ne 0$ for $1 \leq i \leq n$.
Let $F_i$ be a $v+\alpha$-twisted semi-stable sheaf
with $v(F_i)=\pm w(v_i)$ according as the sign of $\rk (w(v_i))$.
If $\rk (w(v_i))>0$, then
\begin{equation}
C_i=\{E\in M_H^{v+\alpha}(v)|\Hom(F_i,E) \ne 0 \}.
\end{equation}
If $\rk (w(v_i))<0$, then
\begin{equation}
C_i=\{E\in M_H^{v+\alpha}(v)|\Hom(E,F_i) \ne 0 \}.
\end{equation}
\end{prop}

\subsection{Normalness of $\overline{M}_H^v(v)$}
\label{subsect:normal}

\begin{prop}\label{prop:normal}
$\overline{M}_H^v(v)$ is normal.
\end{prop}

\begin{proof}
%
We take $\alpha \in D$ with $|\langle \alpha^2 \rangle| \ll 1$.
Let $T$ be a smooth curve and we consider a flat family
of polarized K3 surfaces $\pi:({\cal X},{\cal H}) \to T$
such that 
\begin{enumerate}
\item
$({\cal X}_{t_0},{\cal H}_{t_0})=(X,H)$, $t_0 \in T$,
\item there are families of Mukai vectors 
${\bf v} \in R^* \pi_* {\Bbb Z}$, ${\bf a} \in R^* \pi_* {\Bbb Q}$
 with ${\bf v}_{t_0}=v$,  
${\bf a}_{t_0}={\alpha}$
and 
\item
$\rk \Pic({\cal X}_t) \leq 3$ for a point $t \in T$,
\end{enumerate}
where $({\cal X}_{t},{\cal H}_{t}):=
({\cal X} \otimes k(t),{\cal H} \otimes k(t))$ and
$k(t)$ is the residue field at $t \in T$.
Replacing $T$ by a suitable covering of $T$,
we may assume that there is a section of $\pi$ and 
a locally free sheaf ${\cal E}$ on ${\cal X}$
with $v({\cal E}_t)={\bf v}_t$, $t \in T$.
We consider the relative quot-scheme
$Q:=\Quot_{{\cal E}(-n{\cal H})^{\oplus N}/{\cal X}/T}^{\bf v} \to T$ 
parametrizing all quotients
${\cal E}_t(-n{\cal H}_t)^{\oplus N} \to F$, $t \in T$ with
$v(F)={\bf v}_t$, where $N:=\chi(F(n{\cal H}_t))$.
We denote the universal quotient sheaf by
${\cal F}$. 
We set 
\begin{equation}
Q^{ss}:=\{q \in Q|\text{ ${\cal F}_q:={\cal F} \otimes k(q)$ is
${\bf v}_t$-twisted semi-stable with respect to
${\cal H}_t$ } \}.
\end{equation}
For $n \gg 0$, we have a relative coarse moduli space
$\overline{M}_{{\cal X}/T,{\cal H}}^{\bf v}({\bf v}):=Q^{ss}/PGL(N) \to T$.
Since $T$ is defined over a field of characteristic 0,
$\overline{M}_{{\cal X}/T,{\cal H}}^{\bf v}({\bf v})_t=
\overline{M}_{{\cal H}_t}^{{\bf v}_t}({\bf v}_t)$
(cf. \cite[Thm. 1.1]{MFK:1}).
We also have a relative moduli space 
$\overline{M}_{{\cal X}/T,{\cal H}}^{{\bf v}+{\bf a}}({\bf v}) \to T$.
Replacing $T$ by an open subscheme,
we may assume that 
$\overline{M}_{{\cal X}/T,{\cal H}}^{{\bf v}+{\bf a}}({\bf v})_t$
consists of ${\bf v}_t+{\bf a}_t$-twisted stable sheaves on ${\cal X}_t$
for all $t \in T$ and
there is no walls between ${\bf v}_t$ and
${\bf v}_t+{\bf a}_t$. 
Then $\overline{M}_{{\cal X}/T,{\cal H}}^{{\bf v}+{\bf a}}({\bf v}) \to T$
is a smooth morphism (\cite[Thm. 1.17]{Mu:3}) and  
we have a morphism 
$\Phi:\overline{M}_{{\cal X}/T,{\cal H}}^{\bf v+a}({\bf v}) \to 
\overline{M}_{{\cal X}/T,{\cal H}}^{\bf v}({\bf v})$.
\begin{claim}\label{claim:1}
$\overline{M}_{{\cal X}/T,{\cal H}}^{\bf v}({\bf v})$ is normal.
\end{claim}
Proof of Claim \ref{claim:1}:
It is sufficient to show that $Q^{ss}$ is normal.
By Serre's criterion, we shall show that $Q^{ss}$ is Cohen-Macaulay and
$Q^{ss}$ is regular in codimension 1.
We first prove that $Q^{ss}$ is Cohen-Macaulay.
Let $Q^{spl}$ be the open subscheme of $Q^{ss}$ parametrizing simple sheaves:
\begin{equation}
 Q^{spl}:=\{q \in Q^{ss}| \text{ ${\cal F}_q$ is simple }\}.
\end{equation}
Then $Q^{spl} \to T$ is a smooth morphism (\cite[Thm. 1.17]{Mu:3}).
By the usual deformation theory of sheaves and Lemma \ref{lem:pss} below,
$Q^{ss}_t$ is a locally complete intersection scheme for all $t$.
In particular $Q^{ss}_t$ is Cohen-Macaulay.
Since $Q^{ss}_t$ is smooth at the generic point, 
$Q^{ss}_t$ is reduced. Let $x_t$ be the local parameter of 
$T$ at $t$. By Lemma \ref{lem:ring} below,
$x_t$ is a regular element, which implies that 
$Q^{ss}$ is flat over $T$.
Then $Q^{ss}$ is also Cohen-Macaulay.

We next show that $Q^{ss}$ is regular in codimension 1.
It is sufficient to show that 
$\dim(Q^{ss} \setminus Q^{spl}) \leq \dim Q^{ss}-2$.
By Lemma \ref{lem:pss},
$\dim(Q^{ss}_t \setminus Q^{spl}_t) \leq \dim(Q^{ss}_t \setminus Q^s_t) \leq
 \dim Q^{ss}_t-1$ for all $t \in T$. 
For a point $t \in T$ with $\rho({\cal X}_t) \leq 3$,
by a direct computation, we see that
$\dim(Q^{ss}_t \setminus Q^{spl}_t) \leq \dim Q^{ss}_t-2$.
Since $Q^{spl}$ is an open subscheme of $Q^{ss}$, we get that
$\dim(Q^{ss} \setminus Q^{spl}) \leq \dim Q^{ss}-2$. 
Therefore our claim holds.

Since $\Phi$ is a birational morphism, we get
\begin{equation}\label{eq:Phi}
\Phi_*({\cal O}_{\overline{M}_{{\cal X}/T,{\cal H}}^{\bf v+a}({\bf v})})
={\cal O}_{\overline{M}_{{\cal X}/T,{\cal H}}^{\bf v}({\bf v})}.
\end{equation}
Let $\varphi:\overline{M}_H^{v}(v)_{nor} \to \overline{M}_H^{v}(v)$
be the normalization of $\overline{M}_H^{v}(v)$
and $\phi_{\alpha}':M_H^{v+\alpha}(v) 
\to \overline{M}_H^{v}(v)_{nor}$
the morphism with $\varphi \circ \phi_{\alpha}'=\phi_{\alpha}$.
Since $\overline{M}_H^v(v)_{nor}$ has at worst rational
double points as its singularities, we get $R^1 \phi_{\alpha *}' 
{\cal O}_{M_H^{v+\alpha}(v)}=0$. Since $\varphi$ is a finite morphism,
by the Leray spectral sequence, we get
\begin{equation}
R^1 \phi_{\alpha *} 
{\cal O}_{M_H^{v+\alpha}(v)}=
\varphi_*(R^1 \phi_{\alpha *}' {\cal O}_{M_H^{v+\alpha}(v)})=0.
\end{equation}
Since $\overline{M}_{{\cal X}/T,{\cal H}}^{\bf v+a}({\bf v})$ and
$\overline{M}_{{\cal X}/T,{\cal H}}^{\bf v}({\bf v})$ are flat over $T$,
by using the base change theorem, we get
\begin{equation}
\Phi_*({\cal O}_{\overline{M}_{{\cal X}/T,{\cal H}}^{\bf v+a}({\bf v})}) 
\otimes k(t_0)=
\phi_{{\bf a}_{t_0} *}({\cal O}_{\overline{M}_{{\cal H}_{t_0}}
^{{\bf v}_{t_0}+{\bf a}_{t_0}}({\bf v}_{t_0})}).
\end{equation}
Combining this with \eqref{eq:Phi}, we get 
$\phi_{\alpha *}({\cal O}_{M_H^{v+\alpha}(v)})=
{\cal O}_{\overline{M}_H^v(v)}$, which implies that 
$\overline{M}_H^v(v)$ is normal.
\end{proof}

\begin{lem}\label{lem:pss}
We set 
\begin{equation}
Q^s_t:=\{q \in Q^{ss}_t|\text{ ${\cal F}_q$ is 
${\bf v}_t+{\bf a}_t$-twisted semi-stable} \}.
\end{equation}
Then  
$\dim_{Q^{ss}_t} (Q^{ss}_t \setminus Q^s_t)=\dim Q^s_t-1$.
\end{lem}

\begin{proof}
Let $0 \subset F_1 \subset F_2 \subset \dots \subset
F_s={\cal F}_q$ be the Harder-Narasimhan filtration of 
${\cal F}_q$ with respect to the ${\bf v}_t+{\bf a}_t$-twisted 
semi-stability. 
We set $v_i:=v(F_i/F_{i-1})$.
Then $\langle \widehat{H},v_i \rangle=\langle v,v_i \rangle=0$ and
$\langle v_i^2 \rangle =-2(\rk v_i)^2$.
We shall compute the dimension of
an open subscheme of the flag-scheme 
$F(v_1,v_2,\dots,v_s)$ over $T$
parametrizing filtrations
$0 \subset F_1 \subset F_2 \subset \dots \subset
F_s={\cal F}_q$, $q \in Q_t^{ss}$
such that $F_i/F_{i-1}$, $1 \leq i \leq s$ are ${\bf v}_t+{\bf a}_t$-twisted
semi-stable sheaves with $v(F_i/F_{i-1})=v_i$.
By the arguments in \cite[sect. 3.3]{Y:10}, we get
\begin{equation}
 \dim F(v_1,v_2,\dots,v_s)-\dim GL(N)=\sum_{i>j}\langle v_i,v_j \rangle
 +\sum_i (-(\rk v_i)^2).
\end{equation} 
By the equality
\begin{equation}
2\sum_{i>j}\langle v_i,v_j \rangle
+2\sum_i (-(\rk v_i)^2)=\sum_{i \ne j}\langle v_i,v_j \rangle+
\sum_i \langle v_i^2 \rangle=
\langle v^2 \rangle=0,
\end{equation}
we get 
$\dim F(v_1,v_2,\dots,v_s)-\dim PGL(N)=1$.
Hence our claim holds.
\end{proof}

\begin{lem}\label{lem:ring}
Let $(A,{\frak m})$ and $(B,{\frak n})$ be Noetherian local
rings and $f:A \to B$ a local homomorphism.
Let $x \in {\frak m}$ be a non-zero divisor of $A$ 
satisfying
\begin{enumerate}
\item
$xB={\frak p}_1 \cap {\frak p}_2 \cap \dots \cap
{\frak p}_n$ for some prime divisors ${\frak p}_1,{\frak p}_2,\dots,
 {\frak p}_n$ of $B$, that is, $B/xB$ is reduced,
\item
$B_{{\frak p}_i}$, $1 \leq i \leq n$ are flat over $A$.
\end{enumerate}
Then $x$ is also a non-zero divisor of $B$.
\end{lem}

\begin{proof}
We set
\begin{equation}
K:=\{a \in B| \text{ $x^n a=0$ for some positive integer $n$ } \}.
\end{equation}
$K$ is an ideal of $B$ and
\begin{equation}\label{eq:saturated}
\{a \in B| xa \in K \}=K.
\end{equation}
We shall prove that $K=0$.
By (ii), $K_{{\frak p}_i}=0$ for all $i$.
Since $(K+xB)/xB$ is a sub $B$-module of $B/xB$ and
$B/xB$ is reduced, we get that $(K+xB)/xB=0$.
By \eqref{eq:saturated}, $K=xK$.
By Nakayama's lemma, we get $K=0$.
\end{proof}

\section{Examples}\label{sect:example}

In this section, we shall give some examples of 
$\overline{M}_H^v(v)$ with one singular point.
Let $L:=(-E_8)^{\oplus 2} \oplus U^{\oplus 3}$ be the K3 lattice,
where $U$ is the hyperbolic lattice.
%
\begin{lem}\label{lem:surjective}
Let $N$ be an even lattice of signature $(1,s)$ which has 
a primitive embedding $N \hookrightarrow L$. 
We set $\Delta(N):=\{C \in N|(C^2)=-2\}$.
Assume that there is a primitive element $H$ such that
$(H^2)>0$ and $(H,C) \ne 0$ for all $C \in \Delta(N)$.
Then there is a K3 surface $X$ and an isometry
$f:L \to H^2(X,{\Bbb Z})$ such that
$f(N)=\Pic(X)$ and $f(H)$ is ample.
\end{lem}

\begin{proof}
By the surjectivity of the period map, there is a K3 surface $X$ such that
$\Pic(X)=N$.
We set $\Delta(X)^+:=\{C \in \Pic(X)|\text{ $C$ is a $(-2)$-curve} \}$. 
By the Picard-Lefschetz reflections, we can find a Hodge isomotry
$\phi:H^2(X,{\Bbb Z}) \to H^2(X,{\Bbb Z})$ such that 
$(\phi(f(H)),C)>0$ for all $C \in \Delta(X)^+$.
Replacing $f$ by $\phi \circ f$, we can choose $f(H)$ to be  ample.
\end{proof}

\begin{lem}\label{lem:embed}
Let $(a_{i,j})_{i,j=0}^n$ be a Cartan matrix of affine type
$\tilde{A}_n, \tilde{D}_n, \tilde{E}_n$. 
Let $N_1:=((\oplus_{i=0}^n {\Bbb Z} \beta_i) \oplus {\Bbb Z}\sigma, 
(\;\;,\;\;))$ be a lattice such that
$(\sigma^2)=0$, $(\sigma,\beta_0)=1$,
$(\sigma, \beta_i)=0$, $i>0$ and $(\beta_i,\beta_j)=-a_{i,j}$.
Assume that there is a primitive embedding 
$N_1 \hookrightarrow L$.
Then there is a positive integer $d$ and a primitive sublattice 
$N:=(\oplus_{i=0}^n {\Bbb Z} \xi_i, (\;\;,\;\;))$
of $L$ such that $(\xi_i,\xi_j)=-a_{i,j}+2d$.
\end{lem}

\begin{proof}
Since the signature of $N_1$ is $(1,n+1)$,
there is a vector $x \in N_1^{\perp}$ such that $2d:=(x^2)>0$.
We set $\xi_i:=\beta_i+x$. Then $\oplus_{i=0}^n {\Bbb Z}\xi_i$
is a primitive sublattice of $L$ with  
$(\xi_i,\xi_j)=(\beta_i,\beta_j)+(x^2)=-a_{i,j}+2d$.
\end{proof}

\begin{rem}\label{rem:example}
For $\tilde{A}_n, \tilde{D}_n$, $n \leq 18$ or $\tilde{E}_n$, 
there is a primitive sublattice $N_1$ of $L$
(cf. \cite{S-Z:1}).
\end{rem}

Let $C=(a_{i,j})_{i,j=0}^n$ be a Cartan matrix
of affine type $\tilde{A}_n,\tilde{D}_n,\tilde{E}_n$ and
$Q:=(\oplus_{i=0}^n {\Bbb Z} \alpha_i, (\;\;,\;\;))$ the associated root
lattice, that is, $(\alpha_i,\alpha_j)=a_{i,j}$.
Then there is a vector $\delta:=\sum_{i=0}^n a_i \alpha_i$,
$a_i \in {\Bbb Z}$ such that
\begin{equation}\label{eq:imaginary}
 Q^{\perp}:=\{x \in Q|\text{ $(x,y)=0$ for all $y \in Q$ }\}={\Bbb Z}\delta.
\end{equation}
By the classification of the Cartan matrix of affine type,
we may assume that $a_0=1$.

Let $N:=(\oplus_{i=0}^{n} {\Bbb Z}\xi_{i},(\;\;,\;\;))$ be a primitive 
sublattice of $L$ (Lemma \ref{lem:embed}) such that 
\begin{equation}\label{xi}
(\xi_{i},\xi_{j})= -a_{ij}+2ra,
\end{equation}
where $r$ and $a$ are positive integers with $d=ra$. 
We set
\begin{equation}
H:=\sum_{i=0}^{n} a_{i} \xi_{i}.
\end{equation}

\begin{lem}
\begin{enumerate}
\item[(1)]
$(H,\xi_j)=2ra(\sum_{i=0}^n a_i)$ for all $j$.
In particular, $(H^2)=2ra(\sum_{i=0}^n a_i)^2>0$.
\item[(2)]
$H^{\perp}:=\{\xi \in N|(H,\xi)=0 \}$ is negative definite and
$H^{\perp} \cap \Delta(N)=\emptyset$.
\end{enumerate}
\end{lem}

\begin{proof}
By \eqref{eq:imaginary}, $\sum_{i=0}^n a_i a_{i,j}=0$ for all $j$.
Then we see that
\begin{equation}
(H,\xi_{j})=\sum_{i=0}^{n} a_{i} (\xi_{i} ,\xi_{j})
= \sum_{i=0}^{n} a_{i} (-a_{ij}+2ra)
=2ra(\sum_{i=0}^{n} a_{i})
\label{eqn:H-perp}
\end{equation}
for all $j$. Thus the claim (1) holds. 
We next show the claim (2).
By (\ref{eqn:H-perp}), we see that
\begin{equation}
 H^{\perp}=\{\sum_{i=0}^{n} d_{i} \xi_{i}|\; d_i \in {\Bbb Z},
 \sum_{i=0}^{n} d_{i} =0\}=\bigoplus_{i=0}^{n-1} \mathbb{Z}(\xi_{i}-\xi_{i+1}).
\end{equation}
We define a homomorphism
\begin{equation}
\varphi : H^{\perp} \rightarrow Q=\bigoplus_{i=0}^{n} \mathbb{Z} \alpha_{i}
\end{equation}
by sending $\xi_{i}-\xi_{i+1} \in H^{\perp}$ to 
$ \alpha_{i}-\alpha_{i+1} \in Q$.
Obviously $\varphi$ is injective and 
\begin{equation}\label{eq:varphi}
\im \varphi=\bigoplus_{i=0}^{n-1} {\Bbb Z}(\alpha_i-\alpha_{i+1})=
\{\sum_{i=0}^{n} d_{i} \alpha_{i}|\; d_i \in {\Bbb Z},
 \sum_{i=0}^{n} d_{i}=0\}.
\end{equation}
By \eqref{xi},
we see that
\begin{equation}
(\xi_{i}-\xi_{i+1}, \xi_{j}-\xi_{j+1})= 
-(\alpha_{i}-\alpha_{i+1},\alpha_{j}-\alpha_{j+1}).
\end{equation}
Hence $\varphi$ changes the sign of the bilinear forms.
In order to prove our claim, it is sufficient to show the following assertions:
\begin{enumerate}
\item[(a)]
$\im\varphi$ is positive definite. 
\item[(b)]
There is no vector $x \in \im \varphi$ with $(x^2)=2$.
\end{enumerate}

By \eqref{eq:varphi}, $\delta$ does not belong to $\im \varphi$,
which implies that $\im \varphi$ is positive definite.
We next prove the claim (b).
Since $a_{0}=1$, we can take 
$\{\delta,\alpha_{1},\alpha_{2},\dots,\alpha_{n}\}$
as a $\mathbb{Z}$-basis of $Q$.
Let ${\frak g}$ be the finite simple Lie algebra whose root lattice is 
$\oplus_{i=1}^n {\Bbb Z}\alpha_i$.
Assume that an element $x=l\delta+\sum_{i=1}^n m_i \alpha_i$,
$l, m_i \in {\Bbb Z}$ satisfies that $(x^2)=2$.
Then $\sum_{i=1}^n m_i \alpha_i$ becomes a root of 
${\frak g}$. Hence $m_i \geq 0$ for all $i$,
or $m_i \leq 0$ for all $i$.
Since $\theta:=\sum_{i=1}^n a_i \alpha_i$ is the highest root of ${\frak g}$,
$\sum_{i=0}^n a_i>|\sum_{i=1}^n m_i|$, and hence we get
$l(\sum_{i=0}^n a_i)+\sum_{i=1}^n m_i \ne 0$. Thus $x$ does not 
belong to $\im \varphi$.
\end{proof}

Applying Lemma \ref{lem:surjective} to the lattice $N$ (see \eqref{xi}),
we see that there is a polarized K3 surface $(X,H)$ such that
the Picard lattice of $X$ is $N$ with $H=\sum_{i=0}^n a_i \xi$.
We set
\begin{equation}
\begin{split}
v_{i}:=& r+\xi_{i}+a \rho_{X},\; 0 \leq i \leq n,\\
v:=&\sum_{i=0}^{n} a_{i} v_{i}.
\end{split}
\end{equation}
Then we get that
\begin{equation}
\begin{cases} 
\langle v_{i},v_{j}\rangle =-a_{ij}\\
\langle v,v_j \rangle=0\\
\langle \widehat{H},v_j \rangle=0.
\end{cases}
\end{equation}
Let $E_i$ be a $v$-twisted semi-stable sheaf
with $v(E_i)=v_i$ (Proposition \ref{prop:exist}).
If $E_i$ is properly $v$-twisted semi-stable,
then $\rk \Pic(X) > \rk N$, which is a contradiction.
Hence $E_i$ is $v$-twisted stable for all $i$.
Thus $\overline{M}_H^v(v)$ has a rational double point
of type $(a_{i,j})_{i,j=1}^n$.
In particular, Remark \ref{rem:example} implies that
there is a moduli space $\overline{M}_H^v(v)$ which has
a rational double point of type
$A_n,D_n$, $n \leq 18$, or $E_n$.

\section{Appendix}

Finally we treat the wall crossing phenomenon
under a wall $W_u$ with $u \in {\cal U} \setminus {\cal U}'$.
Assume that $\alpha$ belongs to exactly one wall
$W_u$.
Let $F$ be a $v+\alpha$-twisted semi-stable sheaf
with $v(F)=u$.
By our assumption, $F$ must be $v+\alpha$-twisted stable.
We set $\alpha^{\pm}:=\pm \epsilon v+\alpha$.
We consider the Fourier-Mukai transform ${\cal F}_{\cal E}$
in \eqref{eq:kernel}.
Then we see that $\Ext^1(F,E)=0$ for $E \in M_H^{v+\alpha}(v)$:
Indeed, if $\Ext^1(F,E) \ne 0$, then we have a non-trivial extension
$0 \to E \to G \to F \to 0$
and we see that $G$ is a $v+\alpha$-twisted stable.
On the other hand, $\langle v(G)^2 \rangle<-2$, which is a contradiction.
Hence $\Ext^1(F,E) =0$.
Then, we have an isomorphism
\begin{equation}
\begin{matrix}
M_H^{v+\alpha^-}(v) & \to & M_H^{v+\alpha^-}(v')\\
E & \mapsto & {\cal F}_{\cal E}^1(E),
\end{matrix}
\end{equation}
where $v'=R_{v(F)}(v)$.
We also have an isomorphsim
\begin{equation}
\begin{matrix}
M_H^{v+\alpha^+}(v') & \to & M_H^{v+\alpha^+}(v)\\
E & \mapsto & {\cal F}_{\cal E}^1(E). 
\end{matrix}
\end{equation}
Since $\alpha$ belongs to exactly one wall, we get that
$M_H^{v+\alpha^-}(v')=M_H^{v+\alpha}(v')=M_H^{v+\alpha^+}(v')$.
Therefore we get an isomorphism
$M_H^{v+\alpha^-}(v) \to M_H^{v+\alpha^+}(v)$ and under this identification, 
we get $\theta_v^{\alpha^-}=\theta_v^{\alpha^+}$.

\end{document}